**Optimal Dichotomy of Temporal Scales and Boundedness/ Stability of Time – Varying Multidimensional Nonlinear Systems**

Mark A. Pinsky

**Abstract**. This paper develops a new approach to the estimation of the degree of boundedness/stability of multidimensional nonlinear systems with time-dependent nonperiodic coefficients – an essential task in various engineering and natural science applications. Known approaches to assessing the stability of such systems rest on the utility of Lyapunov functions and Lyapunov first approximation methodologies, typically providing conservative and computationally elaborate criteria for multidimensional systems of this category. Adequate criteria of boundedness of solutions to nonhomogeneous systems of this kind are rare in the contemporary literature. Lately, we develop a new approach to these problems which rests on bounding the evolution of the norms of solutions to initial systems by matching solutions of a scalar auxiliary equation we introduced in [1], [2] and [3]. Still, the technique advanced in [3] rests on the assumption that the average of the linear components of the underlying system is defined by a stable matrix of general position. The current paper substantially amplifies the application domain of this approach. It is merely assumed that the time – dependent linear block of the underlying system can be split into slow and fast varying components by application of any smoothing technique. This dichotomy of temporal scales is determined by the optimal criterion reducing the conservatism of our estimates. In turn, we transform the linear subsystem with slow-varying matrix in a diagonally dominant form by successive applications of the Lyapunov transforms. This prompts the development of novel scalar auxiliary equations embracing the estimation of the norms of solutions to our initial systems. Next, we formulate boundedness/ stability criteria and estimate the relevant regions of the underlying systems using analytical and abridged numerical reasoning. Lastly, we authenticate the developed methodology in inclusive simulations.

## 1. Introduction

Assessment of the boundedness and stability of solutions to nonlinear systems with variable and nonperiodic coefficients is a long-standing problem which finds its origin in a diverse arrays of vital applied fields that, for instance, are concerned with the design of observers and controllers. Sufficient conditions of asymptotic stability for such systems were firstly endorsed by A. M. Lyapunov [4] via the utility of the Lyapunov exponents which were widely used subsequently. In fact, it was shown in [4] that under some natural conditions on the underlying homogeneous system, its trivial solution is asymptotically stable if the linearized at zero system is regular and its maximal Lyapunov exponent is negative, see contemporary review and historical perspectives on this subject, for instance, in [5] and [6]. It turns out that the validation of the regularity condition is a challenging problem and the computation of Lyapunov exponents requires hefty numerical simulations for multidimensional systems. Subsequently, Perron [7] conformed that the regularity condition is essential for the asymptotic stability of such nonlinear systems and observed that the Lyapunov exponents of a linearized system can be sensitive to small perturbations.

A more robust but conservative methodology for stability analysis of nonautonomous nonlinear systems was introduced in [8], where the concept of generalized exponents was developed. It turns out that the upper generalized exponent is equal or larger than the maximal Lyapunov exponent, which extends the conservatism of this approach. Yet, the estimation of generalized exponents for nonautonomous and nonlinear systems presents a challenging problem as well. A more relax stability condition resting on the application of this methodology were given in [9].

The concept of exponential dichotomy of solutions to time – varying systems was used in qualitative analysis and characterization of the stability of some of these systems [8]. However, practical applications of such criteria typically require comprehension of the fundamental set of solutions to the relevant linearized system.

In the control literature, stability conditions for some nonlinear systems with variable coefficients were developed through the application of Lyapunov function method [10] – [23]. Yet, the utility of this approach turns out to be more elaborated for nonautonomous systems where adequate time - dependent Lyapunov functions are rarely available, especially for multidimensional systems of practical concern.

---

Mark A. Pinsky, Department of Mathematics and Statistics, University of Nevada.Reno, Reno NV 89557, USA, e-mial:pinsky@unr.edu.

It appears that there are no necessary and sufficient criteria of local stability of the trivial solution to a sufficiently broad class of time – varying nonlinear systems, see, e.g, [5], [6], [8] – [12].

Several numerical techniques for the estimation of the stability regions of autonomous nonlinear systems were developed in the last few decades [25-27], but extension of these techniques to time – varying systems presently remain elusive. To our knowledge, the problem of estimating the trapping regions for to a wide class of nonhomogeneous nonlinear systems with variable coefficients has not been ascribed in the current literature yet.

An alternative approach to the analysis of the stability of a nonlinear system was developed in [28]. This approach rests on the analysis of convergence of nearby trajectories of the initial system. Yet, to our knowledge, the application of this technique to complex nonlinear systems is still limited.

In [1], we developed a novel approach for assessing the boundedness/stability of time-varying nonlinear systems and applied it to estimation of the trapping/stability regions of these systems. This approach rests on the development of a scalar auxiliary equation with solutions bounding from the above the norm of matching solutions to the original system. The obtained estimates appear to be sufficiently accurate if the solutions are emanated from the central part of the trapping/stability regions but turned out to be more conservative near the boundaries of such regions. Subsequently, [2] combines this technique with the relevant successive approximations granting estimation of the boundaries of the trapping/stability regions and the norms of adjacent solutions with increasing precision.

Yet, the auxiliary equation developed in [1] contains the condition number of the transition matrix of a linear counterpart of the initial system. Such functions can approach infinity on large time intervals, which would make our inferences over conservative for a number of practically important systems. [3] escapes this limitation through developing a technique casting the auxiliary equation in a modified form under the assumption that the average of the time – dependent matrix of a linear block of the given system is a stable matrix of general position.

The current paper develops a novel methodology which essentially extends the application domain of the approach developed in [3]. It merely assumes that the linear block of the given system can be split into slow and fast varying components by application of, e.g., the moving averaging technique. Such temporal dichotomy is defined by the optimal criterion reducing the conservatism of our subsequent estimates. In turn, we apply recursively the Lyapunov transform to map a matrix of the slow linear subsystem to its diagonally - dominant form which ultimately aids the development of novel scalar auxiliary equations. Solutions to these scalar equations bound from the above time – histories of the norms of solutions for a broad class of nonautonomous nonlinear systems and unfold a more inclusive assessment of their boundedness/stability properties. Lastly, we authenticate our inferences in representative simulations.

## 2. Definitions and Preliminaries

This paper studies the behavior of solutions to a time – varying nonlinear system with marked linearization which can be written in the following form,

$$\begin{aligned} &\dot{x} = B(t)x + f_*(t,x) + F_*(t), \quad \forall t \ge t_0, \quad x(t,t_0,x_0) \in \mathbb{R}^n, \quad f_*(t,0) = 0 \\ &x(t_0,t_0,x_0) = x_0, \; x_0 \in H \subset \mathbb{R}^n \end{aligned} \tag{2.1}$$

where functions , $F_*:[t_0,\infty) \to \mathbb{R}^n$, and matrix $B:[t_0,\infty) \to \mathbb{R}^{n\times n}$ are continuous, $H \in \mathbb{R}^n$ is neighborhood of zero, $t_0 \in \mathrm{T} := [\zeta,\infty)$, $\zeta \in \mathbb{R}$, $F_*(t) = F_0\eta(t)$ ,$\sup_{t\ge t_0}\|\eta(t)\| = 1$, $F_0 \ge 0$, $\|\cdot\|$ stands for induced 2-norm of a matrix or 2-norm of a vector, and $x(t,t_0,x_0):[t_0,\infty)\times\mathrm{T}\times H \to \mathbb{R}^n$ is a solution to (2.1). We will frequently write below that $x(t,t_0,x_0) \equiv x(t,x_0)$ to simplify notation and assume that (1.1) possesses a unique solution $\forall x_0 \in H$ and $\forall t \ge t_0$.

The last assumption together with continuity of the right-side of (2.1) implies that $x(t,x_0)$ is a continuous and continuously differential function that is bounded $\forall t \in [t_0,t_*],\; t_* < \infty$ which let us to focus this paper on behavior of $x(t,x_0)$ for $t \to \infty$.

We also write a homogeneous counterpart to (2.1) and its linear counterpart as,

$$\dot{x} = B(t)x + f_*(t,x),\ \forall t \ge t_0$$
$$x(t_0, x_0) = x_0 \in H \subset \mathbb{R}^n \qquad (2.2)$$

,

$$\dot{x} = B(t)x$$
$$x(t_0, x_0) = x_0 \in H \subset \mathbb{R}^n \qquad (2.3)$$

At this point, we assume that equations (2.2) and (2.3) possess unique solutions for $\forall x_0 \in H$ and $\forall t \ge t_0$.

Clearly, the solution to (2.3) takes the following form, $x(t,x_0) = W(t,t_0)x_0$, where $W(t,t_0) = w(t)w^{-1}(t_0)$ and $w(t)$ are the transition and fundamental matrices for (2.3), respectively, which, due to continuity of $B(t)$, are continuous and continuously differentiable in $t$.

Next, we mention that the Lyapunov exponents of solutions to (2.2), which are widely used to assess the exponential growth/decay of solutions to differential equations, can be defined as follows [5], [8],

$$\phi(x_0) = \lim_{t\to\infty}\sup(1/t)\ln\|x(t,x_0)\|$$

After that, we acknowledge the definition of the Lyapunov transform, $v:[t_0,\infty) \to \mathbb{R}^n$, where $v(t)$ is a continuous and continuously differentiable function with $\sup_{t\ge t_0}\left(v(t) + v^{-1}(t) + \dot{v}(t)\right) < \infty$. The Lyapunov transform, $x = v(t)y$ maps (2.3) into a linear system, $\dot{y} = \left(v^{-1}(t)B(t)v(t) - v^{-1}(t)\dot{v}(t)\right)y$ which can possess some desirable properties.

In order, let us bring the definition of the comparison principle [11] which is repeatedly used below. Consider a scalar differential equation,

$$\dot{u}_1 = g(t,u_1),\ \forall t \ge t_0,\ u_1(t,u_{1,0}) \in \mathbb{R}$$
$$u_1(t_0,u_{1,0}) = u_{1,0}$$

where function $g(t,u_1)$ is continuous in $t$ and locally Lipschitz in $u_1$ for $\forall u_1 \in \wp \subset \mathbb{R}$. Suppose that the solution to this equation, $u_1(t,u_{1,0}) \in \wp,\ \forall t \ge t_0$. Next, consider a differential inequality,

$$D^+u_2 \le g(t,u_2),\ \forall t \ge t_0$$
$$u_2(t_0,u_{2,0}) = u_{2,0} \le u_{1,0}$$

where $D^+u_2$ denotes the upper right-hand derivative in $t$ of $u_2(t,u_{2,0})$ and $u_2(t,u_{2,0}) \in \wp,\ \forall t \ge t_0$. Then, $u_1(t,u_{1,0}) \ge u_2(t,u_{2,0}),\ \forall t \ge t_0$.

Our subsequent inferences relay on application of the techniques that were developed in our former papers [1] and [3], which are abbreviated in this section to make this paper more inclusive

In [1] under normalizing condition, $\|w(t_0)\| = 1$ we developed a scalar equation,

$$\dot{X} = p(t)X + c(t)\left\|f_*(t,x(t,x_0)) + F_*(t)\right\|$$
$$X(t_0,X_0) = X_0 = \left\|w^{-1}(t_0)x_0\right\| \qquad , \qquad (2.4)$$

(where

$$p(t) = d\left(\ln\|w(t)\|\right)/dt \qquad (2.5)$$

$$c(t) = \|w(t)\|\left\|w^{-1}(t)\right\| = \sigma_{\max}(w)/\sigma_{\min}(w), \qquad (2.6)$$

$X:[t_0,\infty)\to\mathbb{R}_{\geq 0}$ ($\mathbb{R}_{\geq 0}$ is a set of nonnegative real numbers), and show that the norm of the solution to (2.1) is bounded from above by the solution of (2.4) with matching initial conditions, i.e., $\|x(t,x_0)\|\leq X(t,X_0),\ X_0=\|w^{-1}(t_0)x_0\|,\ \forall t\geq t_0$.

Note that $p(t)$ is a continuous function and $c(t)$ is a continuously differentiable function in $t$ due to continuity of $B(t)$.

To represent (2.4) in a tractable form, we also developed in [1] a nonlinear extension of Lipschitz continuity condition which was written as follows,

$$\|f_*(t,x)\|\leq L(t,\|x\|),\ \forall x\in\Omega_x\subset\mathbb{R}^n,\ \forall t\geq t_0 \tag{2.7}$$

where $L:[t_0,\infty)\times\mathbb{R}_{\geq 0}\to\mathbb{R}_{\geq 0}$ is continuous function in $t$ and $\|x\|$, $L(t,0)=0$, and $\Omega_x$ is a bounded neighborhood of $x\equiv 0$. Note that $L(t,\|x\|)$ can be readily defined in a closed-form if $f_*$ is either a piece-wise polynomial in $x$ or can be approximated by such function with bounded in $\Omega_x$ error term [1]. In the former case, $\Omega_x\equiv\mathbb{R}^n$ and in the latter case, $\Omega_x\equiv\mathbb{R}^n$ if the error term is also globally bounded for $\forall x\in\mathbb{R}^n$. Thus, the condition, $\Omega_x\equiv\mathbb{R}^n$ is readily met in numerous nonlinear systems that are appeared in applications.

Let us show how to define $L(t,\|x\|)$ for a simple system which can be readily generalized to a piecewise polynomial vector – field. Assume that $x=[x_1\ x_2]^T$, $f:\mathbb{R}^2\to\mathbb{R}^2$ and $f=\left[a_1(t)x_1^2x_2^2\quad a_2(t)x_1x_2^2\right]^T$. Then, $\|f\|_2\leq\|f\|_1\leq|a_1||x_1^2||x_2^2|+|a_2||x_1||x_2^2|\leq|a_1|\|x\|^4+|a_2|\|x\|^3$, where we used that $|x_i^n|\leq\|x\|^n$, $i=1,2,\ n\in\mathbb{N}$, $\mathbb{N}$ is a set of positive integers. Note that additional and more complex examples of this kind are provided in Section 7 as well as in [1]-[3].

With (2.7) at hand, (2.4) can be written as follows,

$$\begin{aligned}&\dot{z}=p(t)z+c(t)\left(L(t,z)+\|F_*(t)\|\right)\\&z(t_0,z_0)=\|w^{-1}(t_0)x_0\|=z_0\end{aligned} \tag{2.8}$$

Let us assume that (2.8) possesses a unique solution for $\forall t\geq t_0$ and $\forall z_0\geq 0$ which implies that $z(t,z_0)$ is a continuous and continuously differentiable function $\forall t\geq t_0$ due to the last assumption and continuity of the right side of (2.8).

Clearly, due to comparison principle, $\|x(t,x_0)\|\leq z(t,z_0),\ z_0=\|w^{-1}(t_0)x_0\|,\ \forall t\geq t_0$, where $z(t,z_0)$ is a solution to (2.8). Hence, the last inequality prompts the boundedness/stability criteria for solutions to the multidimensional equations (2.1) and (2.2) through the relevant properties of a scalar equation (2.8).

Yet, these estimates become over conservative for relatively large values of $t$ if $\lim_{t\to\infty}c(t)=\infty$ which frequently takes place even if $A=const$ is a stable matrix.

In [3] we developed a technique lifting this constrain for a practically valuable class of nonlinear systems. Let us review some key points of this technique which will be applied to the current study as well. Firstly, we averaged $B(t)$ as follows, $A(t_0)=\lim_{t_0<t\to\infty}(t-t_0)^{-1}\int_{t_0}^{t}B(s)ds<\infty,\ A(t_0)\in\mathbb{R}^{n\times n},\ \forall t_0\in\mathrm{T}$. Next, we assume that, $A(t_0)\neq 0,\ \|A(t_0)\|<\infty,\ \forall t_0\in\mathrm{T}$, and $A(t_0)$ is a matrix of general position which, hence, is diagonalizable for all but possibly some isolated values of $t_0$. Let us set that complex eigenvalues of $A(t_0)$, $\lambda_k(t_0) = \alpha_k(t_0) \pm i\beta_k(t_0),\ i = \sqrt{-1},\ 1 \le k \le n_1 \le n,\ n_1 \in \mathbb{N}$ and real eigenvalues of $A(t_0)$, $\lambda_k(t_0) = \alpha_k(t_0)$, $n_1 < k \le n$, where $\alpha_k \in \mathbb{R}$ and $\beta_k \in \mathbb{R}_{>0}$, $\mathbb{R}_{>0}$ is a set of positive real numbers. Additionally, we presume that $\alpha_k \ge \alpha_{k+1},\ k \in [1, n-1]$ and also define a square diagonal matrix, $\Lambda = \alpha + i\beta$ with $\alpha = diag(\alpha_1, \alpha_1, ..., \alpha_{n_1}, \alpha_{n_1}, \alpha_{n_1+1}, ..., \alpha_n),\ \beta = diag(\beta_1, -\beta_1, ..., \beta_{n_1}, -\beta_{n_1}, 0, ..., 0)$.

Next, we write (2.1) as follows,

$$\begin{aligned} &\dot{x} = Ax + G_*(t)x + f_*(t,x) + F_*(t), \quad t \in [t_0, \infty) \\ &x(t_0, x_0) = x_0 \end{aligned} \tag{2.9}$$

where $G_*(t) = B(t) - A$ is a zero-mean matrix and, subsequently, rewrite the last equation in the eigenbasis of matrix $A$ as,

$$\begin{aligned} &\dot{y} = (\alpha + i\beta)y + G(t)y + f(t,y) + F(t), \quad t \in [t_0, \infty), \quad y \in \mathbb{R}^n \\ &y(t_0, y_0) = y_0 = v^{-1}x_0 \end{aligned} \tag{2.10}$$

where $y = v^{-1}x$, $y \in \mathbb{C}^n$, $\mathbb{C}^n$ is n-dimensional space of complex values, $v \in \mathbb{C}^{n\times n}$ is the eigenmatrix of $A$, $G = v^{-1}G_*v$, $f(t,y) = v^{-1}f_*(t, vy)$, $F(t) = v^{-1}F_*$.

Afterwards, we rearranged the last equation as follows,

$$\begin{aligned} &\dot{y} = (\eta I + i\beta)y + (\alpha - \eta I)y + G(t)y + f(t,y) + F(t),\ \eta \in \mathbb{R} \\ &y(t_0, y_0) = v^{-1}x_0 \end{aligned}, \tag{2.11}$$

and use in (2.5) and (2.6) that $w(t) = \exp(\eta I + i\beta(t))(t - t_0)$, which returns that $p(w(t)) = \eta$ and $c(w(t)) = 1$. Then, the auxiliary equation (2.8) for (2.11) can be written as,

$$\begin{aligned} &\dot{z} = (\eta + \|\alpha - \eta I\|)z + \|G(t)\|z + L(t,z) + \|F(t)\|,\ \forall \eta \in \mathbb{R} \\ &z(t_0, z_0) = z_0 = \|v^{-1}x_0\| \end{aligned} \tag{2.12}$$

where $\|y(t, y_0)\| \le z(t, z_0),\ \forall t \ge t_0, y \in \Omega_y \subset \mathbb{C}^n$, $\Omega_y$ is a neighborhood of $y \equiv 0$ and function $L(t,z): [t_0, \infty) \times \mathbb{R}_{\ge 0} \to \mathbb{R}_{\ge 0}$ is developed through application of the analog of inequality (2.7) to function, $f(t,y) = v^{-1}f_*(t, vy)$, which returns that

$$\|f(t,y)\| \le L(t, \|y\|),\ \forall y \in \Omega_y,\ \forall t \ge t_0, \tag{2.13}$$

see an illustrative example above in this section as well as more inclusive examples in Section 7. Note that for piecewise polynomial vector fields (2.13) holds for $\forall y \subset \mathbb{C}^n$.

In the sequel, $\eta$ was selected through application of the condition, $\min_\eta(\eta + \|\alpha - \eta I\|)$ yielding that $\eta = \alpha_n$, $(\eta + \|\alpha - \eta I\|)_{|\eta = \alpha_n} = \alpha_1$, which brings (2.12) into the following form,

$$\begin{aligned} &\dot{z} = (\alpha_1 + \|G(t)\|)z + L(t,z) + \|F(t)\|,\ \forall t \ge t_0 \\ &z(t_0, z_0) = z_0 = \|v^{-1}x_0\| \end{aligned} \tag{2.14}$$

Lastly, we present a less conservative counterpart to the prior equation as follows,

$$\dot{z} = \left(\alpha_1 + \left\|\hat{G}(t)\right\|\right) z + L(t,z) + \left\|F(t)\right\|, \ \forall t \ge t_0$$
$$z(t_0, z_0) = z_0 = \left\|v^{-1} x_0\right\| \tag{2.15}$$

where $\hat{G} = G - iE$ , $E(t) = \mathrm{Im}\left(diagG(t)\right)$ , $E \in \mathbb{R}^{n\times n}$ .

Afterwards, we proved in [3] that,

$$\left\|x(t,x_0)\right\| \le \left\|v\right\| z(t,z_0), \ z_0 = \left\|v^{-1}x_0\right\|, \ x_0 \in H, \ \forall t \ge t_0 \tag{2.16}$$

where $x(t,t_0)$ is a solution to (2.1) or (2.2) and $z(t,z_0)$ is the matching solution to the scalar equations (2.14) or (2.15).

Clearly, this approach fails if $A(t_0) = 0$ and becomes overconservative if $A(t_0)$ is an unstable matrix. Moreover, the conservatism of the above technique increases under assumption that $\alpha_1 < 0$ , but $|\alpha_1|$ is sufficiently small and $\sup\limits_{\forall t \ge t_0} \left\|\hat{G}(t)\right\|$ is relatively large.

Next, we present for convenience some conventional definitions of trapping/stability regions, see also [3], where such definitions were also adopted.

**Definition 1**. A connected and compact set of all initial vectors, $\Im_1(t_0)$ is called a trapping region of equation (1.1) if condition $x_0 \in \Im_1(t_0)$ implies that $x(t,x_0) \in \Im_1(t_0), \ \forall t > t_0$ .

Clearly, this definition acknowledges that $\Im_1$ is the invariant set of (1.1).

**Definition 2**. A connected and open set of all initial vectors, $\Im_2(t_0)$ , that includes zero-vector, is called a region of stability of the trivial solution to (1.2) if condition $x_0 \in \Im_2(t_0)$ implies that $x(t,x_0)$ is stable.

**Definition 3**. A connected and open set of all initial vectors, $x_0 \in \Im_3(t_0)$ , that includes zero-vector, is called a region of asymptotic stability of the trivial solution to (1.2) if condition, $x_0 \in \Im_3(t_0)$ implies that $\lim\limits_{t\to\infty} x(t,x_0) = 0$ .

## 3. Generalized Auxiliary Equation

This section develops a novel framework deriving a sequence of scalar auxiliary equations for our initial equations (2.1) or (2.2). Our approach splits matrix $B(t)$ into slow and fast varying components via the application of moving averages to the entries of this matrix. Next, we apply successively the Lyapunov transforms to (2.1) to map a subsystem defined by slow-varying matrix of this equation to its simplified form. This aids the development of more general auxiliary equations escaping the limitations of our prior technique [3].

Our approach is sequenced below in a few consecutive steps where our assumptions are outlined in an ad hoc form. Then, we present a statement encapsulating our hypotheses and outcomes.

3.1. Let us split matrix $B(t)$ into slow and fast- varying components, $B(t) = a_0(t,\delta) + g_0(t,\delta)$, where a slow component, $a_0(t,\delta)$ can be seized, for instance, by application to every entry of $B(t)$ of some moving averages with fixed windows, $\delta = \{\delta_{ij}\}, \ \delta \in \mathbb{R}^{n\times n}_{\ge 0}$ . The components of matrix $\delta$ are defined latter in this section via the application of an optimal criterion. Clearly, the fast matrix $g_0(t,\delta) = B(t) - a_0(t,\delta)$ . Note that the applications of different smoothing techniques aiding selection of slow varying components, like, e.g., wavelets, etc., do not alter the underlying structure of our methodology.

Below we frequently adopt shorter notations, i.e., $a_0(t,\delta)=a_0(t)$ and $g_0(t,\delta)=g_0(t)$ if dependence upon $\delta$ of the relevant matrices is inessential.

Next, we assume that a slow matrix $a_0(t)$ is a matrix of general position [29] which is continuously differentiable for an appropriate number of times that is detailed subsequently. Such matrices possess simple eigenvalues and are diagonalizable for all but possibly some isolated values of $t$ which we call irregular in this paper. Complementary values of $t$ correspond to simple eigenvalues of $a_0(t)$ and are called regular. Thus, the eigenvalues and eigenvectors of $a_0(t)$ are continuously differentiable for a relevant number of times at regular values of $t$, see, e.g., [30], [31].

Revision of our technique for irregular values of $t$ is undertaken in Section 6, whereas other sections of this paper assume that all values of $t \ge t_0$ are regular.

Under assumption that matrix $a_0(t)$ is continuously differentiable, the Lyapunov transform $x=v_1(t)y_1(t)$ maps (2.1) into the following equation,

$$\dot{y}_1=\left(a_1(t)+g_1(t)\right)y_1+f_1(t,y_1)+F_1(t),\ \forall t\ge t_0$$
$$y_1(t_0)=v_1^{-1}(t_0)x_0$$

where, $a_1(t)=\Lambda_1(t)+v_1^{-1}(t)\dot{v}_1(t)$, $\Lambda_1=diag\left(\lambda_{01}(t),...,\lambda_{0n}(t)\right)$, $\lambda_{0i}(t),\ i=1,...,n$ are running eigenvalues of matrix $a_0(t)$ and $g_1(t)=v_1^{-1}(t)g_0(t)v_1(t)$, $f_1(t,y_1)=v_1^{-1}(t)f_*\left(t,v_1(t)y_1\right)$, $F_1(t)=v_1^{-1}(t)F_*(t)$ and $\dot{v}_1(t)$ is a continuous matrix due to our assumption on $a_0(t)$. Hence, $\Lambda_1(t)$, $v_1(t)$ and $\dot{v}_1(t)$ can be considered as slow matrices, which implies that $a_1(t)$ is a slow and diagonally dominant matrix since $\sup_{t_0\le t}\|\dot{v}_1(t)\|$ is a relatively small value.

To develop a consecutive approximation, we assume that matrix $a_0(t)$ is twice continuously differentiable, a matrix $a_1(t)$ possesses simple eigenvalues for $\forall t\ge t_0$ and, hence, is diagonalizable by application of Lyapunov transform $y_1(t)=v_2(t)y_2(t)$, where $v_2(t)$ is an eigenmatrix of $a_1(t)$. This conveys the following equation,

$$\dot{y}_2=\left(a_2(t)+g_2(t)\right)y_2+f_2(t,y_2)+F_2(t),\ \forall t\ge t_0$$
$$y_2(t_0)=V_2^{-1}(t_0)x_0$$

where $a_2(t)=\Lambda_2(t)+v_2^{-1}(t)\dot{v}_2(t),\ \Lambda_2=diag\left(\lambda_{11}(t),...,\lambda_{1n}(t)\right)$, $\lambda_{1i}(t),\ i=1,...,n$ are running eigenvalues of matrix $a_1(t)$, $g_2(t)=v_2^{-1}(t)g_1(t)v_2(t)$, $f_2(t,y_2)=v_2^{-1}(t)f_1\left(t,v_2(t)y_2\right)$, $F_2(t)=v_2^{-1}(t)F_1(t),\ V_2=v_1(t)v_2(t)$, and $\dot{v}_2(t)$ is a continuous matrix with a relatively small $\sup_{t_0\le t}\|\dot{v}_2(t)\|$.

In the sequel, we assume that matrix $a_0(t)$ is $k$-times continuously differentiable, matrix $a_{k-1}(t)$ possesses simple eigenvalues for $\forall t\ge t_0$, and $v_k(t)$ is eigenvector matrix of $a_{k-1}(t)$. Then, application of Lyapunov transform $y_{k-1}(t)=v_k(t)y_k(t),\ k\ge 2$ brings the following equation,

$$\dot{y}_k=\left(a_k(t)+g_k(t)\right)y_k+f_k(t,y_k)+F_k(t),\ \forall t\ge t_0,\ k\ge 2$$
$$y_k(t_0)=V_k^{-1}(t_0)x_0$$

where $a_k(t)=\Lambda_k(t)+v_k^{-1}(t)\dot{v}_k(t),\ \Lambda_k=diag(\lambda_{k-1,1}(t),...,\lambda_{k-1,n}(t)),\ \lambda_{k-1,i}(t),\ i=1,...,n$ are running eigenvalues of matrix $a_{k-1}(t)$ and $g_k(t)=v_k^{-1}(t)g_{k-1}(t)v_k(t)$,

$f_k(t,y_k)=v_k^{-1}(t)f_{k-1}(t,v_k(t)y_k)$, $F_k(t)=v_k^{-1}(t)F_{k-1}(t)$, $V_k(t)=\prod_{j=1}^{k}v_j(t)$ and $\dot{v}_k(t)$ is a continuous matrix.

Lastly, we rewrite the last equation as follows,

$$\begin{aligned}\dot{y}_k&=(\Lambda_k(t)+G_k(t))y_k+f_k(t,y_k)+F_k(t),\ \forall t\ge t_0,\ k\ge 1\\ y_k(t_0)&=V_k^{-1}(t_0)x_0\end{aligned}\tag{3.1}$$

where $G_k(t)=v_k^{-1}(t)\dot{v}_k(t)+g_k(t)$.

3.2. At this point, let us develop a scalar auxiliary equation for equation (3.1) by extending a technique that was developed in [3] and abbreviated in the previous section.

Firstly, we write (3.1) alike (2.11) as follows,

$$\begin{aligned}\dot{y}_k&=(\eta_{k-1}(t)I+i\beta_{k-1}(t))y+(\alpha_{k-1}(t)-\eta_{k-1}(t)I)y_k+G_k(t)y_k+f_k(t,y_k)+F_k(t),\ k\ge 1\\ y_k(t_0)&=V_k^{-1}(t_0)x_0\end{aligned}$$

where $\alpha_{k-1}(t)=\mathrm{Re}(\Lambda_k(t)),\ \beta_{k-1}(t)=\mathrm{Im}(\Lambda_k(t)),\ \alpha_{k-1}\in\mathbb{R}^{n\times n},\ \beta_{k-1}\in\mathbb{R}_{\ge 0}^{n\times n}$ , and a continuous function $\eta_{k-1}:[t_0,\infty)\to\mathbb{R}$ will be defined subsequently. Next, we apply to the last equation a technique which derives the auxiliary equation for (2.11). For this sake, we set that a diagonal matrix $w_{k-1}(t)=\exp\left(\int_{t_0}^{t}(\eta_{k-1}(t)I+i\beta_{k-1}(s))ds\right)$. Then, substitution of $w_{k-1}(t)$ into (2.5) and (2.6) yields that $p(t)=\eta_{k-1}(t),\ c(t)\equiv 1$.

Then, as prior, we select $\eta_{k-1}(t)$ through application of the following condition $\min_{\eta_{k-1}}(\eta_{k-1}(t)+\|\alpha_{k-1}(t)-\eta_{k-1}(t)I\|)$, which yields that $\eta_{k-1}(t)=\alpha_{k-1,\max}(t)$, where $\alpha_{k-1,\max}(t)=\max_{i=1,...,n}\mathrm{Re}(\lambda_{k-1,i}(t))$ is a running maximal real part of eigenvalues of matrix $a_{k-1}(t)$.

Lastly, this comprises the following scalar auxiliary equation,

$$\begin{aligned}\dot{z}_k&=(\alpha_{k-1,\max}(t)+\|G_k(t)\|)z_k+L_k(t,z_k)+\|F_k(t)\|,\ \forall t\ge t_0,\ k\ge 1\\ z_k(t_0,z_{k,0})&=z_{k,0}=\|V_k^{-1}(t_0)x_0\|\end{aligned}\tag{3.2}$$

where functions $L_k:[t_0,\infty)\times\mathbb{R}_{\ge 0}\to\mathbb{R}_{\ge 0}$ are continuous in both variables.

Note that $L_k$ is defined through application of inequality (2.13) to a continuous function $f_k(t,y_k)$ which yields that,

$$\|f_k(t,y_k)\|\le L_k(t,\|y_k\|),\ \forall y_k\in\Omega_{y_k}\subset\mathbb{C}^n,\ \forall t\ge t_0,\ k\ge 1$$

where, as prior, $L_k(t,\|y_k\|)$ can be defined in close-form $\forall y_k\in\mathbb{C}^n$ if , e.g., $f_k(t,y_k)$ are piece-wise polynomials in $y_k$ , see illustrative example in Section 2 and a more inclusive example in Section 7 as well as [1] and [3].

Next, as earlier, we write a more efficient counterpart to equation (2.15) as follows,

$$\dot{z}_k = \left(\alpha_{k,\max}(t) + \left\|\bar{G}_k(t)\right\|\right) z_k + L_k(t, z_k) + \left\|F_k(t)\right\|, \ \forall t \ge t_0, \ k \ge 1$$
$$z_k(t_0, z_{k,0}) = z_{k,0} = \left\|V_k^{-1}(t_0)\, x_0\right\| \tag{3.3}$$

where $\bar{G}_k = G_k - iE_k$, $E_k(t) = \operatorname{Im}\left(diagG_k(t)\right)$, $E_k \in \mathbb{R}^{n\times n}$.

To simplify further referencing, we present a homogeneous counterpart to (3.3) as follows,

$$\dot{z}_k = \left(\alpha_{k,\max}(t) + \left\|\bar{G}_k(t)\right\|\right) z_k + L_k(t, z_k), \ \forall t \ge t_0, \ k \ge 1$$
$$z_k(t_0, z_{k,0}) = z_{k,0} = \left\|V_k^{-1}(t_0)\, x_0\right\| \tag{3.4}$$

3.3. Subsequently, we recall that all entries of the right side of equations (3.2) – (3.4) depend upon entries of matrix $\delta$ controlling separation of temporal scales of matrix $B(t)$. While the above equations are defined for $\forall \delta \in \mathbb{R}_{\ge 0}^{n\times n}$, the efficacy of the estimates they provide can be enhanced by appropriate choice of $\delta$. To this end, to abate the influence of nonlinear component in the right sides of our auxiliary equations in a sufficiently small neighborhood of $z_k \equiv 0$, we assume that $L_k(t, z_k, \delta) \le \varsigma_{k,1}(t,\delta)\, z_k^{\varsigma_{k,2}}$, where continuous $\varsigma_{k,1} : [t_0, \infty) \times \mathbb{R}_{\ge 0}^{n\times n} \to \mathbb{R}_{\ge 0}$ and scalars $\varsigma_{k,2} > 1$. Next, we define $\delta$ to maximize the degrees of stability of either equation, $\dot{z}_K = \left(\alpha_{K,\max}(t,\delta) + \left\|G_K(t,\delta)\right\|\right) z_K$ or equation, $\dot{z}_K = \left(\alpha_{K,\max}(t,\delta) + \left\|\bar{G}_K(t,\delta)\right\|\right) z_K$ for some $k = K \ge 1$.

For this sake, we define Lyapunov exponents of solutions to these equations as follows,

$$\phi_{K,1}(\delta) = \limsup_{T\to\infty} \left( T^{-1} \int_0^T \left(\alpha_{K,\max}(t,\delta) + \left\|G_K(t,\delta)\right\|\right) dt \right)$$
$$\phi_{K,2}(\delta) = \limsup_{T\to\infty} \left( T^{-1} \int_0^T \left(\alpha_{K,\max}(t,\delta) + \left\|\bar{G}_K(t,\delta)\right\|\right) dt \right)$$

and select $\delta$ to conform one of the following relations $\inf_{0\le\delta\le\Delta} \phi_{K,1}(\delta)$ or $\inf_{0\le\delta\le\Delta} \phi_{K,2}(\delta)$, where entries of matrix $\Delta \subset \mathbb{R}_{\ge 0}^{n\times n}$ assign some practically feasible intervals of variation of components of matrix $\delta$. The efficacy of this approach is assessed in the simulations presented in Section 7.

3.4. Lastly, we encapsulate our inferences in the following**,**

**Theorem 1**. Suppose that continuous matrix $B(t) = a_0(t) + g_0(t)$, where $a_0(t)$ and $g_0(t)$ are slow and fast varying matrices, respectively. Assume also that matrix $a_0(t)$ is $k$ - times continuously differentiable, matrices $a_k(t) = \Lambda_k(t) + v_k^{-1}(t)\dot{v}_k(t), \ k \ge 1, \ \forall t \ge t_0$ possess simple eigenvalues $\forall t \ge t_0$, functions $L_k(t, z_k), \ k \ge 1$ are continuous in both variables $\forall z_k > 0$ and $\forall t \ge t_0$, function $F_*(t)$ is continuous $\forall t \ge t_0$, and $\Omega_{y_k} \equiv \mathbb{C}^n$. Also assume that equations (2.1), (2.2) and (3.2) – (3.4) possess unique solutions $\forall t \ge t_0$, $\forall x_0 \in H$, $\forall z_{k,0} \ge 0$, respectively.

Then, the right sides of equations (3.2) – (3.4) are continuous functions in the relevant variables and the following inequality holds,

$$\left\|x(t, x_0)\right\| \le \left\|V_k(t)\right\| z_k(t, z_{k,0}), \ z_{k,0} = \left\|V_k^{-1}(t_0)\, x_0\right\|, \quad \forall x_0 \in H, \ \forall t \ge t_0, \ k \ge 1 \tag{3.5}$$

where $x(t, x_0)$ is a solution to either (2.1) or (2.2) and $z_k(t, z_{k,0})$ is a solution to one of the corresponding scalar equations (3.2), (3.3), or (3.4), respectively.

**Proof**. The proof of this statement directly follows from our previous derivation. In fact, the implied conditions assure that matrices $a_k(t),\ k \ge 1$ are diagonalizable and their eigenvalues and eigenvectors are continuously differentiable [30], [31]. This implies that $\alpha_k(t)$ and, hence, $\alpha_{k,\max}(t)$, and other terms in the right side of (3.2) – (3.4) are at least continuous functions.

Next, diagonalization of matrices $a_k(t)$ by the appropriate Lyapunov transforms leads to (3.1) and consequently forges the design of auxiliary equations (3.2) or (3.3) $\forall y_k \in \mathbb{C}^n$, which prompts inequality (3.5)□

Apparently, there is no warranty that inequality (3.5) becomes sharper for every successive value of $k$ which reminisced alike behavior of asymptotic series. Still, to prove the boundedness/stability of solutions to our initial multidimensional equations, one should attest a required property for at least one scalar auxiliary equation with $k = K$.

Our approach is essentially simplified if $B(t)$ is already a slow matrix. In this case, it let to conform and extend some criteria of stability of the trivial solution to equation (2.2) which were developed via the utility of nonlinear versions of so known frozen coefficients approach, see, e.g., [22], [23]. Moreover, our methodology also embraces the boundedness criteria for nonhomogeneous systems of this kind which practically were overlooked in the current literature.

Let us note as well that equations (3.2) and (3.3) are more general than our prior auxiliary equations (2.14) and (2.15) and approach their former counterparts as $\delta_{ij} \to \infty$.

Local analysis of the boundedness and stability of solutions to equations (3.3) and (3.4) is abridged and can be carried out in nearly closed form under the following conditions,

$$L_k(t, z_k) \le l_k(t, \hat{z}_k) z_k,\ \ \forall t \ge t_0,\ \ z_k \in [0, \hat{z}_k],\ \hat{z}_k > 0,\ k \ge 1 \tag{3.6}$$

where functions $l_k; [t_0, \infty) \times \mathbb{R}_{\ge 0} \to \mathbb{R}_{\ge 0}$ are continuous in $t$ and bounded in $\hat{z}_k$. Clearly, (3.6) can be interpreted, for instance, as an application of Lipschitz continuity condition to a scalar nonnegative function $L_k(t, z_k)$. More efficient definition of functions $l_k(t, \hat{z}_k)$ can be readily given if $L_k(t, z_k)$ are piece – wise polynomial in $z_k$, see, e.g., [3] for further details.

Note that alike notations used in this paper and in [3] let us in the next two sections naturally extend the inferences that we developed prior to more general types of nonlinear and nonautonomous systems.

## 4. Local Stability and Boundedness Criteria

This section utilizes (6.3) to derive simplified stability and boundedness criteria and abridge the estimations of boundedness/stability regions.

Under presumption that all values of $t \ge t_0$ are regular, application of (3.6) to (3.3) yields a sequence of scalar linear equations which can be written as follows,

$$\begin{aligned} &\dot{Z}_k = \mu_k(t, \hat{z}_k) Z_k + \left\| F_k(t) \right\|,\ \ Z_k \in [0, \hat{z}_k],\ \forall t \ge t_0,\ k \ge 1 \\ &Z_k(t_0, z_{k,0}, \hat{z}_k) = z_{k,0} = \left\| V_k^{-1}(t_0) x_0 \right\| \end{aligned} \tag{4.1}$$

where functions $\mu_k(t, \hat{z}_k) = \alpha_{k,\max}(t) + \left\| \bar{G}_k(t) \right\| + l_k(t, \hat{z}_k)$ are continuous in $t$ and bounded in $\hat{z}_k$, and functions $\left\| F_k(t) \right\| = F_0 \left\| V_k^{-1}(t) \eta(t) \right\|$ are continuous.

The solutions to (4.1) take the form,

$$Z_k(t, z_{k,0}, \hat{z}_k) = z_{k.0} Z_{k,h}(t, \hat{z}_k) + F_0 Z_{k,F}(t, \hat{z}_k),\ k \ge 1 \tag{4.2}$$

where $Z_{k,h}(t,t_0,\hat{z}_k)=\exp\left(\int_{t_0}^{t}\mu_k(\tau,\hat{z}_k)d\tau\right)$, $Z_{k,F}(t,t_0,\hat{z})=\int_{t_0}^{t}\theta_k(t,\tau,\hat{z})\left\|V_k^{-1}(\tau)\eta(\tau)\right\|d\tau$ and

$\theta_k(t,\tau,\hat{z}_k)=\exp\left(\int_{\tau}^{t}\mu_k(s,\hat{z}_k)ds\right)$. Note that below we frequently adopt a reduced notation, i.e.,

$Z_{k,h}(t,t_0,\hat{z}_k)=Z_{k,h}(t,\hat{z}_k)$ and $Z_{k,F}(t,t_0,\hat{z})=Z_{k,F}(t,\hat{z})$.

Clearly, $Z_{k,h}(t,\hat{z}_k)>0,\ \forall t\ge t_0,\ \forall\hat{z}_k>0$ and $\theta_k(t,\tau,\hat{z})>0,\ \forall t>\tau,\ \forall t\ge t_0,\ \forall\hat{z}_k>0$, and the last inequality yields that $Z_{k,F}(t,\hat{z})>0,\ \forall t\ge t_0,\ \forall\hat{z}_k>0$ .

### 4.1. Stability Criteria for Time – Varying Homogeneous Nonlinear Systems

Next, let us write a homogeneous counterpart to (4.1)

$$\begin{aligned}&\dot{Z}_k=\mu_k(t,\hat{z}_k)Z_k,\ Z_k\in[0,\hat{z}_k],\ \forall t\ge t_0,\ k\ge 1\\&Z_k(t_0,z_{k,0},\hat{z}_k)=z_{k,0}=\left\|V_k^{-1}(t_0)x_0\right\|\end{aligned}\tag{4.3}$$

with the following solutions, $Z_k(t,z_{k,0},\hat{z}_k)=z_{k.0}Z_{k,h}(t,\hat{z}_k)$.

As is known, the stability/asymptotic stability of equation (4.3) can be described through conditions set on the fundamental solution matrices for these equations, see, e.g., [24] and additional references therein. Application of these conditions yields that equation (4.3) with $k=K\ge 1$ is stable if and only if,

$$\int_{t_0}^{t}\mu_K(\tau,\hat{z}_K)d\tau<\infty,\ \forall t\ge t_0,\ \forall\hat{z}\in(0,\hat{z}_{K,b}),\ \forall t_0\in\mathrm{T}\tag{4.4}$$

and is asymptotically stable if and only if

$$\lim_{t\to\infty}\sup\left(\int_{t_0}^{t}\mu_K(\tau,\hat{z}_K)d\tau\right)=-\infty,\ \forall\hat{z}\in(0,\hat{z}_{K,B}),\ \forall t_0\in\mathrm{T}\tag{4.5}$$

where $\hat{z}_{K,b}$ and $\hat{z}_{K,B}$ are the values of $\hat{z}_K$ for which either (4.4) or (4.5) hold.

Next lemma provides sufficient conditions prompting the linearization of (3.4) by (3.6).

**Lemma.** Assume that $\exists k=K\ge 1$ such that a function $\mu_K(t,\hat{z}_K)$ is continuous in $t$ and bounded in $\hat{z}_k$ $\forall t\ge t_0,\ \forall\hat{z}_k>0$ and (4.3) is stable $\forall\hat{z}_K\in(0,\hat{z}_{K,\max})$, where $\hat{z}_{K,\max}$ is the maximal value of $\hat{z}_K$ for which (4.3) is stable that can be infinity. Also assume that $z_{K,0}<\hat{z}_K/Z_{K,h,s}(t_0,\hat{z}_K)$, where $Z_{K,h,s}(t_0,\hat{z}_K)=\sup_{t\ge t_0}Z_{K,h}(t,t_0,\hat{z}_K)$.

Then,

$$\sup_{t\ge t_0}Z_K(t,z_{K,0},\hat{z}_K)<\hat{z}_K,\ \forall t\ge t_0,\ \forall z_{K,0}<\hat{z}_K/Z_{K,h,s}(t_0,\hat{z}_K),\ \hat{z}_K\in(0,\hat{z}_{K,\max})\tag{4.6}$$

**Proof.** Evidently, $\sup_{t\ge t_0}Z_{K,h}(t,t_0,\hat{z}_K)=Z_{K,h,s}(t_0,\hat{z}_k)<\infty$ since $Z_{K,h}(t,t_0,\hat{z})$ is continuous in $t$ and bounded in $\hat{z}_k$, and (4.3) is stable, which yields that $\sup_{t\ge t_0}Z_K(t,t_0,z_{K,0},\hat{z}_K)=Z_{K,h,s}(t_0,\hat{z}_K)z_{K,0}<\infty$ and, consequently, (4.6) holds □

Note that condition $z_{K,0} < \hat{z}_K / Z_{K,h,s}(t_0, \hat{z}_K)$ is automatically assured if $z_{K,0}$ is sufficiently small, $\hat{z}_K > 0$ and (4.3) is stable.

**Theorem 2.** Assume that $\exists k = K \geq 1$ such that assumptions of Theorem 1 and the Lemma are met, and equation (4.3) is stable/ asymptotically stable $\forall \hat{z}_K \in (0, \hat{z}_{K,\max})$. Then, the trivial solution to (2.2) is stable/ asymptotically stable, respectively, and in both cases inequality (3.5) takes the form,

$$\begin{aligned} &\|x(t, x_0)\| \leq \|V_K(t)\| z_K(t, z_{K,0}) \leq \|V_K(t)\| Z_K(t, z_{K,0}, \hat{z}_K), \ \forall t \geq t_0, \\ &\forall z_{K,0} < \hat{z}_K / Z_{K,h,s}(t_0, \hat{z}_K), \ x_0 \mid \|V_K^{-1}(t_0) x_0\| < \hat{z}_K / Z_{K,h,s}(t_0, \hat{z}_K), \ \hat{z}_K \in (0, \hat{z}_{K,\max}) \end{aligned} \tag{4.7}$$

where $Z_K(t, z_{K,0}, \hat{z}_K)$, $z_K(t, z_{K,0})$ and $x(t, x_0)$ are solutions to equations (4.3), (3.4) and (2.2), respectively.

Furthermore, the region of stability or asymptotic stability of the trivial solution to (2.2) includes the following set initial vectors,

$$\|V_K^{-1}(t_0) x_0\| < \sup_{\hat{z}_K \in (0, \hat{z}_{K,\max})} \left( \hat{z}_K / Z_{K,h,s}(t_0, \hat{z}_K) \right) \tag{4.8}$$

Clearly, (4.8) defines the interior of the relevant ellipsoid

**Proof**. Firstly, let us show that $z_K(t, z_{K,0}) < \hat{z}_K, \ \forall t \geq t_0, \ \hat{z}_K \in (0, \hat{z}_{K,\max})$ if $z_{K,0} < \hat{z}_K / Z_{K,h,s}(t_0, \hat{z}_K)$ and (4.3) is stable. Pretend in contrary that $t_* > t_0$ is the smallest value of $t$ such that $z_K(t_*, z_{K,0}) = \hat{z}_K$ under prior conditions. Then, (3.6) and, thus, equation (4.3) holds for $\forall t \in [t_0, t_*]$ which, due to comparison principle and (4.6), implies that $z_K(t_*, z_{K,0}) \leq Z(t_*, z_{K,0}, \hat{z}_K) < \hat{z}_K$ . Thus, $z_K(t, z_{K,0}) < \hat{z}_K, \ \forall t \geq t_0, \ \forall z_{K,0} < \hat{z}_K / Z_{K,s}$ which enables application of the linear equation (4.3). Then, the application of comparison principal prompts that $z_K(t, z_{K,0}) \leq Z(t, z_{K,0}, \hat{z}_K), \ \forall t \geq t_0, \ z_{K,0} < \hat{z}_K / Z_{K,h,s}, \ \hat{z}_K \in (0, \hat{z}_{K,\max})$ which implies that the trivial solution to (3.4) is stable.

Due to the last inequality, (3.5) can be written as (4.7) which implies stability of the trivial solution to (2.2) under the conditions of this theorem.

Let us subsequently estimate the regions of stability of the trivial solutions to equations (3.4) and, in turn, (2.2). Clearly, a solution to (3.4) is stable if $z_{K,0} < \sup_{\hat{z}_K \in (0, \hat{z}_{K,\max})} \left( \hat{z}_K / Z_{K,s}(t_0, \hat{z}_K) \right)$ which permits linearization of (3.4). This implies that the region of stability of the trivial solution to (2.2) is defined by (4.8) since $z_{K,0} = \|V_K^{-1}(t_0) x_0\|$.

In the sequel, we assume that (4.3) is asymptotically stable and other conditions of this statement hold. Then, under this more conservative condition, our prior inferences and (4.7) hold, which implies that $\lim_{t \to \infty} Z_K(t, z_{K,0}, \hat{z}_K) = 0$ and, subsequently, $\lim_{t \to \infty} \|x(t, x_0)\| = 0$ □

The above statement prompts some apparent sufficient stability conditions of the trivial solution to (2.2) which can be readily checked.

**Corollary 1**. Assume that $\exists k = K \geq 1$ such that conditions of Theorem 1 and Lemma are met and:

(1) either (4.4) or (4.5) hold or

(2) $$\mu_K(t, \hat{z}_K) \leq -\nu_K(\hat{z}_K), \ \forall t \geq t_0, \ \nu_K(\hat{z}_K) > 0, \ \forall \hat{z}_K \in (0, \hat{z}_{K,v}) \tag{4.9}$$

where $\hat{z}_{K,v}$ is the maximal value of $\hat{z}_K$ for which (4.9) holds. Then, (1) the trivial solution to (2.2) is either stable or asymptotically stable, respectively; and (2) the trivial solution of (2.2) is asymptotically stable. Furthermore, in both cases the inequalities (4.7) and (4.8) hold, where $\hat{z}_{K,\max}$ should be changed on either $\hat{z}_{K,b}$ or $\hat{z}_{K,B}$ or $\hat{z}_{K,v}$ in (4.7) and (4.8), respectively.

**Proof**. Really, under the above conditions, (4.3) is either stable or asymptotically stable, which, due to Theorem 2, assures this statement □

Evidently, the Lyapunov exponents for solutions to equation (4.3) can be defined as follows,

$\phi_k(\hat{z}_k) = \limsup_{t\to\infty}(t-t_0)^{-1}\int_{t_0}^{t}\mu_k(\tau,\hat{z}_k)\,d\tau$, which prompts,

**Corollary 2**. Assume that $\exists k = K \ge 1$ such that conditions of Theorem 1 and Lemma are met, and $\phi_K(\hat{z}_K) < 0,\ \hat{z}_K \in (0, \hat{z}_{K,\phi})$, where $\hat{z}_{K,\phi}$ is the maximal value of $\hat{z}_K$ for which last inequality holds.

Then, the trivial solution to (2.2) is asymptotically stable and the inequalities (4.7) and (4.8) hold, where $\hat{z}_{K,\max}$ should be exchanged on $\hat{z}_{K,\phi}$ in (4.7) and (4.8).

**Proof**. In fact, the condition of this statement assumes that the Lyapunov exponent of solutions to equation (4.3) with $k = K$ is negative which yields that the trivial solution to this equation is asymptotically stable [5] and, due to Theorem 2, (4.7) and (4.8) hold □

### 4.2 Boundedness Criteria for Solutions of Nonhomogeneous Time-Varying Nonlinear Systems

Now we are going to develop some boundedness criteria of solutions to (2.1) and estimate the relevant regions of the initial data. These conditions frequently assure either the stability or asymptotic stability of the forced solutions to (2.1) and, thus, are important in applications. Our inferences have some correspondence with the concept of input-to-state stability developed by E. D. Sontag [33], [34] under more restrictive conditions.

A solution to (4.1) $Z_K(t, t_0, z_{K,0}, \hat{z}_K) < \infty,\ \forall t \ge t_0,\ \hat{z}_K \in (0, \hat{z}_{K,1})$ if the trivial solution to (4.3) is either stable or asymptotically stable and

$$Z_{K,F,s}(t_0, \hat{z}_K) = \sup_{t\ge t_0} Z_{K,F}(t, t_0, \hat{z}_K) < \infty,\ \hat{z}_K \in (0, \hat{z}_{K,F}),\ \forall t_0 \in \mathrm{T}, \tag{4.10}$$

where $\hat{z}_{K,F}$ is the maximal value of $\hat{z}_K$ for which (4.10) holds, and $\hat{z}_{K,1} = \min(\hat{z}_{K,\max}, \hat{z}_{K,F})$.

This leads us to

**Theorem 3**. Assume that $\exists k = K \ge 1$ such that assumptions of Theorem 1 hold, equation (4.3) is stable $\forall \hat{z}_K \in (0, \hat{z}_{K,\max})$, function $\mu_K(t, \hat{z}_K)$ is continuous in $t$ and bounded in $\hat{z}_K$, and (4.10) and inequality

$$z_{K,0} < \left(\hat{z}_K - F_0 F_{K,F,s}(t_0, \hat{z}_K)\right) / Z_{K,h,s}(t_0, \hat{z}_K) = z_{K,F}(t_0, \hat{z}_K),\ \hat{z}_K \in (0, \hat{z}_{K,1}) \tag{4.11}$$

hold. Then,

(I) $\|x(t, x_0)\| \le \infty, \forall t \ge t_0$, where $x(t, x_0)$ are solutions to (2.1) stemming from the interior of the ellipsoid in $\mathbb{R}^n$ which is defined as follows,

$$x_0 \mid \|V_K^{-1}(t_0)x_0\| < \sup_{\hat{z}_K \in (0, \hat{z}_{K,1})} z_{K,F}(t_0, \hat{z}_K) \tag{4.12}$$

(II) inequality (3.5) takes the form,

$$\begin{gathered}\|x(t, x_0)\| \le \|V_K(t)\| Z_K(t, z_{K,0}, \hat{z}_K) \le \|V_K(t)\|\left(z_{K,0} Z_{K,h,s}(t_0, \hat{z}_K) + F_0 Z_{K,F,s}(t_0, \hat{z}_K)\right),\\ \forall t \ge t_0,\ x_0 \mid \|V_K^{-1}(t_0)x_0\| < z_{K,F},\ z_{K,0} < z_{K,F},\ \hat{z}_K \in (0, \hat{z}_{K,1})\end{gathered} \tag{4.13}$$

where $Z(t, z_{K,0}, \hat{z}_K)$ and $x(t, x_0)$ are solutions to (4.1) and (2.1), respectively, and

(III) if $\exists k = K$ such that (4.3) is asymptotically stable, other conditions of this statement hold, and $\limsup_{t\to\infty}\|V_K(t)\| = \bar{V}_K < \infty$; then,

$$\limsup_{t\to\infty}\|x(t, x_0)\| \le F_0 \bar{V}_K Z_{K,F,s}(t_0, \hat{z}_K),\ x_0 \mid \|V_K^{-1}(t_0)x_0\| < z_{K,F},\ \hat{z}_K \in (0, \hat{z}_{K,1}) \tag{4.14}$$

**Proof**. Firstly, we show that under the condition of this theorem,

$$\sup_{t\ge t_0} Z_K\left(t, z_{K,0}, \hat{z}_K\right) < z_{K,0} Z_{K,h,s}\left(t_0, \hat{z}_K\right) + F_0 Z_{K,F,s}\left(t_0, \hat{z}_K\right) < \hat{z}_K,$$
$$\forall t \ge t_0,\ \forall z_{K,0} < z_{K,F} \tag{4.15}$$

As is mentioned prior, $Z_{K,h,s}\left(t_0, \hat{z}_K\right) < \infty$ since we assume that for $k = K$ (4.3) is stable. Then,

$$\sup_{t\ge t_0} Z_K\left(t, t_0, z_{K,0}, \hat{z}_K\right) \le z_{K,0} \sup_{t\ge t_0} Z_{K,h}\left(t, t_0, \hat{z}_K\right) + \sup_{t\ge t_0} Z_{K,F}\left(t, t_0, \hat{z}_K\right) =$$
$$z_{K,0} Z_{K,h,s}\left(t_0, \hat{z}_K\right) + F_0 Z_{K,F,s}\left(t_0, \hat{z}_K\right)$$

since $Z_{K,h}\left(t, t_0, \hat{z}_K\right) \ge 0$ and $Z_{K,F}\left(t, t_0, \hat{z}_K\right) \ge 0$ which implies (4.15).

Next, we show that $z_K\left(t, z_{K,0}\right) < \hat{z}_K,\ \forall t \ge t_0,\ \forall \hat{z}_K \in \left(0, \hat{z}_{K,1}\right)$ if (4.3) is stable and $z_{K,0} < z_{K,F}$. In fact, pretend, in contrary, that $t_* > t_0$ is the smallest value of $t$ such that $z_K\left(t_*, z_{K,0}\right) = \hat{z}_K$. Then, (3.6) and, thus, equation (4.1) holds for $\forall t \in \left[t_0, t_*\right]$. Hence, the comparison principle and (4.15) yield that $z_K\left(t_*, z_{K,0}\right) \le Z\left(t_*, z_{K,0}, \hat{z}_K\right) < \hat{z}_K$. This contradiction authenticates that $z_K\left(t, z_{k,0}\right) < \hat{z}_K,\ \forall t \ge t_0,\ \forall z_{K,0} < z_{K,F}$, $\hat{z}_K \in \left(0, \hat{z}_{K,1}\right)$ which enables application of equation (4.1) in our inferences.

Consequently, application of the comparison principal returns that,
$z_K\left(t, z_{K,0}\right) \le Z_K\left(t, z_{K,0}, \hat{z}_K\right) < \infty,\ \forall t \ge t_0,\ z_{K,0} < z_{K,F},\ \hat{z}_K \in \left(0, \hat{z}_{K,1}\right)$.

Clearly, the last inequality let us rewrite (3.5) as (4.13), which shows that the norms of solutions to (2.1) are bounded if they are outgoing from the region specified by (4.12).

Lastly, (4.14) directly follows from (4.13) if equation (4.3) is asymptotically stable for $k = K$ since in this case $\lim_{t\to\infty} Z_{K,h,s}\left(t, x_0\right) = 0$ □

Next, two corollaries present simplified sufficient conditions under which the presumptions of Theorem 4 are met.

**Corollary 3**. Assume that $\exists k = K$ such that conditions of Theorem 1are met, the inequalities (4.9), (4.10), and (4.11) hold, and function $\mu_K\left(t, \hat{z}_K\right)$ is continuous in $t$ and bounded in $\hat{z}_K$. Then, solutions to (2.1) are bounded in the region that is defined by (4.12) and inequalities (4.13) and (4.14) hold if $\hat{z}_{K,1} = \min\left(\hat{z}_{K,v}, \hat{z}_{K,F}\right)$
.

**Proof**. In fact, (4.9) implies asymptotic stability of (4.3) for $\forall \hat{z}_K \in \left(0, \hat{z}_{K,v}\right]$ which, together with other assumptions, assure that the conditions of theorem 3 hold □

In the sequel, let us assume that the Lyapunov exponent of solutions to (4.3), i.e., $\phi_K\left(\hat{z}_K\right) < 0,\ \hat{z}_K \in \left(0, \hat{z}_{K,\phi}\right)$, see Corollary 2 for more details. This implies that $Z_{K,h}\left(t, z_{K,0}, \hat{z}_K\right) \le D_K\left(\varepsilon_K, \hat{z}_K\right) e^{-\chi_K(t-t_0)}$ and that $\theta_K\left(t, \tau, \hat{z}_K\right) \le D_K\left(\varepsilon_K, \hat{z}_K\right) e^{-\chi_K(t-\tau)}$, where $-\chi_K = \phi_K + \varepsilon_K,\ \chi_K > 0$ and $\varepsilon_K > 0$ are small values, and constants $D_K\left(\varepsilon_K, \hat{z}_K\right) > 0$ [8].

This leads to

**Corollary 4**. Assume that $\exists k = K$ such that conditions of theorem 1 are met, the inequalities (4.10) and (4.11) hold, function $\mu_K\left(t, \hat{z}_K\right)$ is continuous in $t$ and bounded in $\hat{z}_K$, and $\phi_K\left(\hat{z}_K\right) < 0,\ \forall \hat{z}_K \in \left(0, \hat{z}_{K,\phi}\right)$. Then, solutions to (2.1) are bounded in the region that is defined by (4.12), and inequalities (4.13) and (4.14) hold for $\hat{z}_{K,1} = \min\left(\hat{z}_{K,\phi}, \hat{z}_{K,F}\right)$.

**Proof**. In fact, in this case (4.3) is asymptotically stable [5], which implies that the conditions of theorem 4 hold □

Note that the above statements grant boundedness/stability conditions of solutions to equations (2.1) or (2.2) in nearly closed form.

**5. Nonlocal Boundedness/Stability Criteria**

Under assumption that that all values of $t \ge t_0$ are regular, this section devises the boundedness or stability criteria of the initial equations (2.1) and (2.2) from the analysis of auxiliary equations (3.3) or (3.4), respectively. Alike notations, used in equations (3.3) and (2.14), foster a formal similarity between our current and prior criteria that were developed in [3]. Yet, the current criteria are substantially broader since they rely on applications of more general auxiliary equation.

Next theorem accentuates that under some natural conditions, the solutions to (3.3) monotonically increase in the initial value for $\forall t \ge t_0$ which abridges simulations of this equation.

**Theorem 4**. Assume that equations (3.3)/ (3.4) admit unique solutions $\forall z_{k,0} \ge 0,\ k \ge 1$. Also assume that $z_k\left(t, z'_{k0}\right)$ and $z_k\left(t, z''_{k0}\right)$ are solutions to (3.3)/(3.4) with $0 \le z'_{k0} \le z''_{k0}$ . Then, $z\left(t, z'_{k0}\right) \le z\left(t, z''_{k0}\right),\ \forall t \ge t_0$.

**Proof**. In fact, this theorem directly follows from the uniqueness property of solutions to the scalar equations (3.3)/(3.4) □

The prior statement leads to the boundedness/stability criteria for our initial equations which can be tested via reduced numerical simulations of (3.3)/ (3.4). To formulate these criteria, let us firstly define a set of centered at zero concentric ellipsoids, $E(z) \subset \mathbb{R}^n$ as follows,

$$E(z): \left\| V^{-1}(t_0)\, x \right\| = z \ge 0 \tag{5.1}$$

Also, we assume that $\partial E(z) \subset \mathbb{R}^{n-1}$ defines the boundaries of these ellipsoids and $E_-(z) = E(z) - \partial E(z)$ . This prompts the following

**Theorem 5**. 1.Assume that $\exists k = K \ge 1$ such that the trivial solution to the relevant equation (3.4) is either stable or asymptotically stable. Then the trivial solution to (2.2) is stable/asymptotically stable as well. Furthermore, let the interval $[0, \bar{z}_K)$ defines the region of stability/asymptotic stability of the trivial solution to (3.4). Then, the set $E_-\left(\bar{z}_K\right)$ is enclosed in the stability region of (2.1).

2. Assume that $\exists k = K \ge 1$ such that the interval $[0, \bar{z}_K]$ defines a trapping region for the corresponding equation (3.3), i.e., $z_K\left(t, z_{K0}\right) \le \bar{z}_K, \forall z_{K0} \in [0, \bar{z}_K],\ \forall t \ge t_0$ . Then, ellipsoid $E\left(\bar{z}_k\right)$ is enclosed into the trapping region of equation (2.1), i.e., $\left\| x\left(t, x_0\right) \right\| < \infty,\ \forall x_0 \in E\left(\bar{z}_K\right),\ \forall t \ge t_0$.

**Proof**. Really, the proof of both statements follows from inequality (3.5), where it is presumed that $x\left(t, x_0\right)$ and $z\left(t, z_{K0}\right)$ are solutions to either equations (2.2) and (3.4) or (2.1) and (3.3), respectively □

Thus, to estimate the trapping/stability regions of multidimensional equations (2.1) or (2.2) about $x \equiv 0$ we, firstly, shall distinguish in simulations the threshold value, $\bar{z}_K$ which defines the relevant intervals of the initial values for equations (3.3) or (3.4). This task is markedly streamlined since $z_K\left(t, z_{K0}\right)$ is monotonically increases in $z_{K0},\ \forall t \ge t_0$ due to theorem 4. Then the ellipsoids estimating the trapping/stability regions for equations (2.1) or (2.2) are defined by the following formula, $\bar{z}_K = \left\| V^{-1}(t_0)\, x_0 \right\|$.

To get further insight into the behavior of solutions to our nonlinear auxiliary equations, we can bound from the above and below the right sides of (3.3) or (3.4) by their time-invariant counterparts which introduce two scalars, autonomous and, thus, integrable equations. Solutions to these equations define bilateral bounds on the relevant solutions of (3.3) or (3.4) which evoke an analytical assessment of the behavior of the solutions to our initial equations. This approach was realized in [3] using more restrained auxiliary equations but can be immediately extended to a wider class of systems corresponding to the auxiliary equations developed in this paper.

## 6. Auxiliary Equation in a Neighborhood of an Irregular Time-moment

Currently let us presume that $a_k(t),\ k \ge 0$ are matrices of general position which are not diagonalizable for some isolated irregular values of $t$ where their normal forms possess two-dimensional Jordan- blocks corresponding to repeated eigenvalues, see [29], pp. 219-235.

To simplify further reasoning, let us pretend that matrices $a_k(t),\ k \ne k_1$ possess simple eigenvalues $\forall t \ge t_0$ but a matrix $a_{k_1}(t)$ is not diagonalizable at $t = \tilde{t}_{k_1}$ and $\Delta_{k_1} = \left[\tilde{t}_{k_1} - \breve{t}_{k_1},\ \ \tilde{t}_{k_1} + \breve{t}_{k_1}\right],\ \breve{t}_{k_1} > 0$ is a small interval about this isolated irregular time moment. Thus, our prior technique fails to provide the auxiliary equation at $t = \tilde{t}_{k_1}$, $\forall k > k_1$.

Nonetheless, under some additional conditions, we can define the auxiliary equation in a neighborhood of an isolated irregular time moment using the mean value theorem for matrices which we adopt in the following form

$$a_{k_1}(t) = a_{k_1}(t+d) - d\int_0^1 a'_{k_1}\left(t + d(1-s)\right)ds,\ t \in \Delta_{k_1} \tag{6.1}$$

where a small real value $d \ne 0$ and $a'_{k_1}(t)$ is the derivative of matrix $a_{k_1}(t)$.

To unfold this procedure, we firstly acknowledge that a matrix $a_{k_1}(t+d),\ \forall t \in \Delta_{k_1},\ d \ne 0$ possesses simple eigenvalues and, alike as it was done in Section 3, successively simplify the underlying equations using the Lyapunov transforms, $y_k(t+d) = v_{k_1,k+1}(t+d)\,y_{k+1}(t+d),\ \ \forall t \in \Delta_{k_1},\ \forall k \ge k_1$. This yields matrices $a_{k_1,k}(t+d) = \Lambda_{k_1,k}(t+d) + v_{k_1,k}^{-1}(t+d)\dot{v}_{k_1,k}(t+d),\ t \in \Delta_{k_1},\ \forall k > k_1$, where $\Lambda_{k_1,k} = diag\left(\lambda_{k_1,k,1}(t+d),...,\lambda_{k_1,k,n}(t+d)\right)$, $\lambda_{k_1,k,j}(t+d),\ j = 1,...,n$ are running eigenvalues of matrices $a_{k_1,k}(t+d)$, and $\alpha_{k_1,k,\max}(t+d)$ is the maximal real part of eigenvalues $\lambda_{k_1,k,j}(t+d)$.

Next, $\forall k > k_1$ and $\forall t \in \Delta_{k_1}$, we set that $V_{k_1,k}(t+d) = \prod_{m=k_1+1}^{k} v_{k_1,m}(t+d)$,

$g_{k_1,k}(t,d) = V_{k_1,k}^{-1}(t+d)\,g_{k_1}(t)\,V_{k_1,k}(t+d)$, $f_{k_1,k}(t,y_k,d) = v_{k_1,k}^{-1}(t+d)\,f_{k_1,k-1}\left(t, v_{k_1,k}(t+d)\,y_k\right)$,

$F_{k_1,k}(t,d) = V_{k_1,k}^{-1}(t+d)\,F_{k_1}(t)$, $G_{k_1,k}(t,d) = v_k^{-1}(t+d)\,\dot{v}_k(t+d) + g_{k_1,k}(t,d)$,

$M_{k_1,k}(t,d) = dV_{k_1,k}^{-1}(t+d)\left(\int_0^1 a'_{k_1}\left(t+d(1-s)\right)ds\right)V_{k_1,k}(t+d),\ \ \forall t \in \Delta_{k_1}$.

Next, as in Section 3, we assume that a function $L_{k_1,k}(t,z_k),\ t \in \Delta_{k_1},\ z_k \ge 0,\ k > k_1$ is defined as in Section 3 to content inequality, $\left\| f_{k_1,k}(t,y_k,d)\right\| \le L_{k_1,k}\left(t,\|y_k\|,d\right),\ \forall y_k \in \mathbb{C}^n$, $\forall t \in \Delta_{k_1}$ and also recall that the Lyapunov transforms introduced in Section 3, i.e., $y_{k-1}(t) = v_k(t)\,y_k(t)$ remain intact in the current section $\forall t \notin \Delta_{k_1}$ and $k \ge 1$ which leads to the following,

**Theorem 6**. Suppose that matrix $a_{k_1}\left(\tilde{t}_{k_1}\right),\ t_{k_1} > t_0$ is not diagonalizable but all other assumptions of Theorem 1 are met and matrix $a_0(t)$ is $k+1$- times continuously differentiable $\forall t \in \Delta_{k_1}$. Also assume that matrices $a_{k_1}(t+d)$ and $a_{k_1,k}(t+d) = \Lambda_{k_1,k}(t+d) + v_{k_1,k}^{-1}(t+d)\,\dot{v}_{k_1,k}(t+d),\ k > k_1$ are matrices of general position which possess simple eigenvalues $\forall t \in \Delta_{k_1}$ and some small $d \ne 0$ and that functions $L_{k_1,k}(t,z_k,d),\ \forall t \in \Delta_{k_1},\ z_k \ge 0,\ k > k_1$ are continuous in $t$, $z_k$ and $d$.

Then,

1. Equation (3.2) remains intact $\forall k \le k_1$ and $\forall t \ge t_0$ .
2. $\forall k > k_1$, (3.2) can be altered as follows,

$$\dot{\underleftarrow{z}}_k = \left(\alpha_{k-1,\max}\left(t\right) + \left\|G_k\left(t\right)\right\|\right)\underleftarrow{z}_k + L_k\left(t, \underleftarrow{z}_k\right) + \left\|F_k\left(t\right)\right\|, \quad \forall t \in \left[t_0, \tilde{t}_{k_1} - \breve{t}_{k_1}\right]$$

$$\underleftarrow{z}_k\left(t_0,\right) = z_{k0}$$

$$\dot{\underleftrightarrow{z}}_k = \left(\alpha_{k_1,k-1,\max}\left(t+d\right) + \left\|G_{k_1,k}\left(t,d\right)\right\| + \left\|M_{k_1,k}\left(t,d\right)\right\|\right)\underleftrightarrow{z}_k +$$

$$L_{k_1,k}\left(t, \underleftrightarrow{z}_k, d\right) + \left\|F_{k_1,k}\left(t,d\right)\right\|, \; \forall t \in \Delta_{k_1} \tag{6.2}$$

$$\underleftrightarrow{z}_k\left(\tilde{t}_{k_1} - \breve{t}_{k_1}\right) = \underleftarrow{z}_k\left(\tilde{t}_{k_1} - \breve{t}_{k_1}\right)$$

$$\dot{\underrightarrow{z}}_k = \left(\alpha_{k-1,\max}\left(t\right) + \left\|G_k\left(t\right)\right\|\right)\underrightarrow{z}_k + L_k\left(t, \underrightarrow{z}_k\right) + \left\|F_k\left(t\right)\right\|, \quad \forall t > \tilde{t}_{k_1} + \breve{t}_{k_1}$$

$$\underrightarrow{z}_k\left(\tilde{t}_{k_1} + \breve{t}_{k_1}\right) = \underleftrightarrow{z}_k\left(\tilde{t}_{k_1} + \breve{t}_{k_1}\right)$$

where $\alpha_{k_1,k-1,\max}$ are continuous functions $\forall t \in \Delta_{k_1}$ and matrices $G_{k_1,k}$ and $M_{k_1,k}$ are continuous and, hence, are bounded in norm $\forall t \in \Delta_{k_1}$ .

3. Assume additionally that the initial value problem (6.2) possesses a unique solution

$$\underline{z}_k\left(t, z_{k,0}\right) = \begin{cases} \underleftarrow{z}_{k,1}\left(t, z_{k,0}\right), & \forall t \in \left[t_0, \tilde{t}_{k_1} - \breve{t}_{k_1}\right] \\ \underleftrightarrow{z}_k\left(t, z_{k,0}\right), & \forall t \in \left[\tilde{t}_{k_1} - \breve{t}_{k_1} \quad \tilde{t}_{k_1} + \breve{t}_{k_1}\right] \\ \underrightarrow{z}_k\left(t, z_{k,0}\right), & \forall t > \tilde{t}_{k_1} + \breve{t}_{k_1} \end{cases}$$

Then, inequality (3.5) takes the following form,

$$\left\|x\left(t, x_0\right)\right\| \le \left\|V_k\left(t\right)\right\| \underline{z}_k\left(t, z_{k,0}\right), \quad \forall x_0 \in H, \; z_{k,0} = \left\|V_k^{-1}\left(t_0\right) x_0\right\|, \; \forall t \ge t_0, \; k > k_1 \tag{6.3}$$

**Proof**. $k+1$-times differentiability of $a_0\left(t\right), \; \forall t \in \Delta_{k_1}$ together with assumptions on the structure of the eigenvalues of the relevant matrices, imply that matrices $v_{k_1,k}\left(t+d\right)$ , $a_{k_1,k}\left(t+d\right)$, $G_{k_1,k}\left(t,d\right)$ and functions $\alpha_{k_1,k-1,\max}\left(t+d\right)$, and $F_{k_1,k}\left(t\right)$ are at least continuous $\forall k > k_1$ $\forall t \in \Delta_{k_1}$ . This assures that $\alpha_{k_1,k-1,\max}\left(t+d\right) < \infty, \; \forall t \in \Delta_{k_1}$ , $\left\|G_{k_1,k}\left(t,d\right)\right\| < \infty, \; \forall t \in \Delta_{k_1}$ and $\left\|F_{k_1,k}\left(t\right)\right\| < \infty, \; \forall t \in \Delta_{k_1}$ . In turn, the same assumptions yields that a matrix $\sigma\left(t\right) = \int_0^1 a'_{k_1}\left(t + d\left(1-s\right)\right) ds$ is continuous and, thus, $\left\|\sigma\left(t\right)\right\| < \infty, \; \forall t \in \Delta_{k_1}$ .

Next, using (6.1), we write the underlying equation for $k = k_1$ and $\forall t \in \Delta_{k_1}$ as follows,

$$\dot{\underleftrightarrow{y}}_{k_1} = \left(a_{k_1}\left(t+d\right) - d\int_0^1 a'_{k_1}\left(t + d\left(1-s\right)\right) ds + g_{k_1}\left(t\right)\right)\underleftrightarrow{y}_{k_1} + f_{k_1}\left(t, \underleftrightarrow{y}_{k_1}\right) + F_{k_1}\left(t\right), \; \forall t \in \Delta_{k_1}$$

$$\underleftrightarrow{y}_{k_1}\left(\tilde{t}_{k_1} - \breve{t}_{k_1}\right) = y_{k_1}\left(\tilde{t}_{k_1} - \breve{t}_{k_1}\right)$$

Then, subsequent applications to the last equation of the transforms $y_k\left(t+d\right) = v_{k+1}\left(t+d\right) y_{k+1}\left(t+d\right)$, $\forall t \in \Delta_{k_1}, \; \forall k \ge k_1$ lead to the following equations,

$$\dot{\underleftrightarrow{y}}_k = \left(a_{k_1,k}\left(t+d\right)+g_{k_1,k}\left(t,d\right)+M_{k_1,k}\left(t,d\right)\right)\underleftrightarrow{y}_k + f_{k_1,k}\left(t,\underleftrightarrow{y}_k,d\right)+F_{k_1,k}\left(t,d\right),\ \forall t\in\Delta_{k_1},\ k>k_1$$

$$\underleftrightarrow{y}_k\left(\tilde{t}_{k_1}-\breve{t}_{k_1}\right)=\underleftarrow{y}_k\left(\tilde{t}_{k_1}-\breve{t}_{k_1}\right)$$

In turn, we alter the previous equation as follows

$$\dot{\underleftrightarrow{y}}_k = \left(\Lambda_{k_1,k}\left(t+d\right)+G_{k_1,k}\left(t,d\right)+M_{k_1,k}\left(t,d\right)\right)\underleftrightarrow{y}_k + f_{k_1,k}\left(t,\underleftrightarrow{y}_k,d\right)+F_{k_1,k}\left(t,d\right),\ \forall t\in\Delta_{k_1},\ k>k_1$$

$$\underleftrightarrow{y}_k\left(\tilde{t}_{k_1}-\breve{t}_{k_1}\right)=\underleftarrow{y}_k\left(\tilde{t}_{k_1}-\breve{t}_{k_1}\right)$$

Consequently, application of the technique used in Section 3 to the last equation, brings the second equation in (6.2) which assures application of inequality (6.3) $\forall t\le\tilde{t}_{k_1}+\breve{t}_{k_1}$.

In turn, inequality $\left\|x\left(t,x_0\right)\right\|\le\left\|V_k\left(t\right)\right\|\underrightarrow{z}_k\left(t,z_{k,0}\right),\ \forall t>\tilde{t}_{k_1}+\breve{t}_{k_1},\ x_0\in H,\ z_{k,0}=\left\|V_k^{-1}\left(t_0\right)x_0\right\|$ holds due to comparison principle, which assures (6.3) $\forall t\ge t_0$ □

Note that immediate extension of the prior reasoning assures Theorem 6 under more relax conditions that some matrices $a_{k_1}\left(t\right),\ \ldots,a_{k_J}\left(t\right),k_1<\ldots<k_J$ are not diagonalizable at the same isolated time moment and also if some of matrices $a_k\left(t\right)$ are not diagonalizable at a finite number of isolated time moments. However, such generalizations of the prior theorem are left beyond the scope of this paper.

Consequently, we extend the application of Theorem 2 to the case concerned in this section. For this sake, we restrain inequality (3.6) on the pertaining time interval as follows,

$$L_{k_1,k}\left(t,z_k,d\right)\le\underleftrightarrow{l}_{k_1,k}\left(t,\hat{z}_k,d\right)z_k,\ \ \forall t\in\Delta_{k_1},\ \ z_k\in[0,\hat{z}_k],\ \ k>k_1 \tag{6.4}$$

where a function $\underleftrightarrow{l}_{k_1,k}\left(t,\hat{z}_k,d\right)$ is continuous in $t$ and bounded in $\hat{z}_k$ and $d$. Clearly, we employed in (3.6) and (6.4) the same value of $\hat{z}_k$ for simplicity which, in some cases, might extend the conservatism of the follow-up inferences.

Applications of (3.6) and (6.4) yields a set of three linear equations with matching initial conditions

$$\begin{aligned}
&\dot{\underleftarrow{Z}}_k=\underleftarrow{\mu}_k\left(t,\hat{z}_k\right)\underleftarrow{Z}_k+\left\|\underleftarrow{F}_k\left(t\right)\right\|,\ \ \forall\underleftarrow{Z}_k\in[0,\hat{z}_k],\ \forall t\in\left[t_0,\tilde{t}_{k_1}-\breve{t}_{k_1}\right],\ \hat{z}_k>0,\ k>k_1\\
&\underleftarrow{Z}_k\left(t_0\right)=z_{k0}\\
&\dot{\underleftrightarrow{Z}}_k=\underleftrightarrow{\mu}_{k_1,k}\left(t,\hat{z}_k,d\right)\underleftrightarrow{Z}_k+\left\|F_{k_1,k}\left(t,d\right)\right\|,\ \ \forall\underleftrightarrow{Z}_k\in[0,\hat{z}_k],\ \forall t\in\Delta_{k_1},\ d\ne 0,\ k>k_1\\
&\underleftrightarrow{Z}_k\left(\tilde{t}_{k_1}-\breve{t}_{k_1}\right)=\underleftarrow{Z}_k\left(\tilde{t}_{k_1}-\breve{t}_{k_1}\right)\\
&\dot{\underrightarrow{Z}}_k=\underrightarrow{\mu}_k\left(t,\hat{z}_k\right)\underrightarrow{Z}_k+\left\|\underrightarrow{F}_k\left(t\right)\right\|,\ \ \forall\underrightarrow{Z}_k\in[0,\hat{z}_k],\ \forall t>\tilde{t}_{k_1}+\breve{t}_{k_1},\ k>k_1\\
&\underrightarrow{Z}_k\left(\tilde{t}_{k_1}+\breve{t}_{k_1}\right)=\underleftrightarrow{Z}_k\left(\tilde{t}_{k_1}+\breve{t}_{k_1}\right)
\end{aligned} \tag{6.5}$$

where functions $\underleftarrow{\mu}_k$, $\underrightarrow{\mu}_k$, $\underleftarrow{F}_k$ and $\underrightarrow{F}_k$ are developed by restraining on the appropriate time intervals of their relevant counterparts that are defined in Section 4, and

$$\underleftrightarrow{\mu}_{k_1,k}\left(t,\hat{z}_k,d\right)=\alpha_{k_1,k-1,\max}\left(t+d\right)+\left\|G_{k_1,k}\left(t,d\right)\right\|+\left\|M_{k_1,k}\left(t,d\right)\right\|+\underleftrightarrow{l}_{k_1,k}\left(t,\hat{z}_k,d\right),\ \forall t\in\Delta_{k_1},\ \forall k>k_1$$

Clearly, $\underleftarrow{z}_k\left(t\right)\le\underleftarrow{Z}_k\left(t\right),\ \forall t>\tilde{t}_{k_1}-\breve{t}_{k_1},\ \forall k>k_1,\ \underleftrightarrow{z}_k\left(t\right)\le\underleftrightarrow{Z}_k\left(t\right),\ \forall t\in\Delta_{k_1},\ \forall k>k_1$ and $\underrightarrow{z}_k\left(t\right)\le\underrightarrow{Z}_k\left(t\right),\ \forall t>\tilde{t}_{k_1}+\breve{t}_{k_1},\ \forall k>k_1$ due to comparison principle.

Next, let us define piecewise continuous functions

$$\underline{\mu}_{k_1,k}\left(t,\hat{z}_k,d\right)=\begin{cases}\underset{\leftarrow}{\mu}_k\left(t,\hat{z}_k\right), & \forall t\in\left[t_0,\tilde{t}_{k_1}-\breve{t}_{k_1}\right],\ \hat{z}_k>0\\ \underset{\leftrightarrow}{\mu}_{k_1,k}\left(t,\hat{z}_k,d\right), & \forall t\in\left[\tilde{t}_{k_1}-\breve{t}_{k_1}\quad \tilde{t}_{k_1}+\breve{t}_{k_1}\right],\ d\neq 0,\ k>k_1\\ \underset{\rightarrow}{\mu}_k\left(t,\hat{z}_k\right), & \forall t>\tilde{t}_{k_1}+\breve{t}_{k_1}\end{cases}$$

and

$$\varphi_{k_1,k}\left(t,d\right)=\begin{cases}\underset{\leftarrow}{F}_k\left(t\right), & \forall t\in\left[t_0,\tilde{t}_{k_1}-\breve{t}_{k_1}\right]\\ F_{k_1,k}\left(t,d\right), & \forall t\in\left[\tilde{t}_{k_1}-\breve{t}_{k_1}\quad \tilde{t}_{k_1}+\breve{t}_{k_1}\right],\ d\neq 0,\ k>k_1\\ \underset{\rightarrow}{F}_k\left(t\right), & \forall t>\tilde{t}_{k_1}+\breve{t}_{k_1}\end{cases}$$

Then, alike (4.2), continuous general solution to (6.5) can be written as follows

$$\underline{Z}_k\left(t,t_0,z_{k,0},\hat{z}_k,d\right)=z_{k.0}\underline{Z}_{k,h}\left(t,t_0,\hat{z}_k,d\right)+F_0 Z_{k,\varphi}\left(t,t_0,\hat{z}_k,d\right),\ k\geq k_1 \tag{6.6}$$

where $\underline{Z}_{k,h}\left(t,t_0,\hat{z}_k,d\right)=\exp\left(\int_{t_0}^{t}\underline{\mu}_{k_1,k}\left(\tau,\hat{z}_k,d\right)d\tau\right)$, $\underline{\theta}_k\left(t,\tau,\hat{z}_k d\right)=\exp\left(\int_{\tau}^{t}\underline{\mu}_{k_1,k}\left(s,\hat{z}_k,d\right)ds\right)$,

$Z_{k,\varphi}\left(t,t_0,\hat{z},d\right)=F_0^{-1}\int_{t_0}^{t}\underline{\theta}_k\left(t,\tau,\hat{z},d\right)\left\|\varphi_{k_1,k}\left(\tau,d\right)\right\|d\tau$ .

Assume subsequently that $\exists k=K\geq k_1$ such that homogeneous counterpart to (6.5) is stable $\forall\hat{z}_K\in\left(0,\hat{z}_{K,M\text{ax}}\right)$. Then,

$$\sup_{t\geq t_0}\underline{Z}_{K,h}\left(t,t_0,\hat{z}_K,d\right)=\underline{Z}_{K,h,s}\left(t_0,\hat{z}_K,d\right)<\infty,\ \forall\hat{z}_K\in\left(0,\hat{z}_{K,M\text{ax}}\right) \tag{6.7}$$

which let us to reframe Theorem 2 as follows.

**Theorem 7**. Assume that the assumptions of Theorem 6 are met and $\exists k=K>k_1$ such that function $l_K\left(t,\hat{z}_K\right)$ is continuous in $t$ and bounded and in $\hat{z}_K$, function $\underset{\leftrightarrow}{l}_{k_1,K}\left(t,\hat{z}_K,d\right)$ is continuous in $t$ and bounded in $\hat{z}_K$ and $d$ , and the homogeneous counterpart of (6.5) is stable or asymptotically stable $\forall\hat{z}_K\in\left(0,\hat{z}_{K,M\text{ax}}\right)$. Then, the trivial solution to equation (2.2) is stable or asymptotically stable, respectively, and in both cases inequality (4.7) takes the form,

$$\begin{aligned}&\left\|x\left(t,x_0\right)\right\|\leq z_{K,0}\left\|V_K\left(t\right)\right\|\underline{Z}_{K,h}\left(t,\hat{z}_k,d\right),\ \forall t\geq t_0,\ \hat{z}_K\in\left(0,\hat{z}_{K,M\text{ax}}\right),\ d>0\\ &\forall z_{K,0}<\hat{z}_K/\underline{Z}_{K,h,s}\left(t_0,\hat{z}_k,d\right),\ \forall x_0\mid\left\|V_K^{-1}\left(t_0\right)x_0\right\|<\hat{z}_K/\underline{Z}_{K,h,s}\left(t_0,\hat{z}_K,d\right)\end{aligned}$$

(6.8)

Furthermore, the region of stability or asymptotic stability of the trivial solution to equation (2.2) includes the following set of initial vectors,

$$\left\|V_K^{-1}\left(t_0\right)x_0\right\|<\sup_{\hat{z}_K\in\left(0,\hat{z}_{K,M\text{ax}}\right)}\left(\hat{z}_K/\underline{Z}_{K,h,s}\left(t_0,\hat{z}_K,d\right)\right) \tag{6.9}$$

Clearly, (6.9) defines the interior region of the relevant ellipsoid

**Proof**. Clearly, (6.6) let to apply directly the steps used to prove Theorem 2 to conform to the current statement □

Adopting Theorem 3 to the assumptions of this section includes the alteration of inequalities (4.10) and (4.11). In fact, let us assume that the homogeneous counterpart of (6.5) is stable or asymptotically stable $\forall\hat{z}_K\in\left(0,\hat{z}_{K,M\text{ax}}\right)$ as well as that,

$$Z_{K,\varphi,s}\left(t_0,\hat{z}_K,d\right)=\sup_{t\ge t_0} Z_{K,\varphi}\left(t,t_0,\hat{z}_K,d\right)<\infty,\ \hat{z}_K\in\left(0,\hat{z}_{K,\varphi}\right),\ \forall t_0\in \mathrm{T} \tag{6.10}$$

and

$$z_{K,0}<\left(\hat{z}_K-F_0F_{K,\varphi,s}\left(t_0,\hat{z}_K,d\right)\right)/\underline{Z}_{K,h,s}\left(t_0,\hat{z}_K,d\right)=\underline{z}_{K,\varphi}\left(t_0,\hat{z}_K,d\right),\ \hat{z}_K\in\left(0,\hat{\underline{z}}_{K,1}\right) \tag{6.11}$$

where $\hat{z}_{K,\varphi}$ is the maximal value of $\hat{z}_K$ for which (6.10) holds and $\hat{\underline{z}}_{K,1}=\min\left(\hat{z}_{K,M\mathrm{ax}},\hat{z}_{K,\varphi}\right)$.

This leads to

**Theorem 8**. Assume that the assumptions of Theorem 6 are met and $\exists k=K>k_1$ such that function $l_K\left(t,\hat{z}_K\right)$ is continuous in $t$ and bounded and in $\hat{z}_K$, function $\underset{\leftrightarrow}{l}_{k_1,K}\left(t,\hat{z}_K,d\right)$ is continuous in $t$ and bounded in $\hat{z}_K$ and $d$, and the homogeneous counterpart to (6.5) is stable or asymptotically stable $\forall \hat{z}_K\in\left(0,\hat{z}_{K,M\mathrm{ax}}\right)$, and inequalities (6.10) and (6.11) hold. Then,

(I) $\left\|x\left(t,x_0\right)\right\|\le\infty,\forall t\ge t_0$, where $x\left(t,x_0\right)$ are solutions to (2.1) stemmed from the interior of the ellipsoid in $\mathbb{R}^n$ which is defined as follows,

$$x_0\mid\left\|V_K^{-1}\left(t_0\right)x_0\right\|<\sup_{\hat{z}_K\in\left(0,\hat{\underline{z}}_{K,1}\right)}\underline{z}_{K,\varphi}\left(t_0,\hat{z}_K,d\right) \tag{6.12}$$

(II) inequality (3.5) takes the form,

$$\begin{aligned}&\left\|x\left(t,x_0\right)\right\|\le\left\|V_K\left(t\right)\right\|\underline{Z}_K\left(t,z_{K,0},\hat{z}_K,d\right)\le\\&\left\|V_K\left(t\right)\right\|\left(z_{K,0}\underline{Z}_{K,h,s}\left(t_0,\hat{z}_K,d\right)+F_0Z_{K,\varphi,s}\left(t_0,\hat{z}_K,d\right)\right),\\&\forall t\ge t_0,\ x_0\mid\left\|V_K^{-1}\left(t_0\right)x_0\right\|<\underline{z}_{K,\varphi}\left(t_0,\hat{z}_K,d\right),\ z_{K,0}<\underline{z}_{K,\varphi},\ \hat{z}_K\in\left(0,\hat{\underline{z}}_{K,1}\right)\end{aligned} \tag{6.13}$$

where $\underline{Z}_K\left(t,z_{K,0},\hat{z}_K,d\right)$ and $x\left(t,x_0\right)$ are solutions to (6.5) and (2.1), respectively, and

(III) if $\exists k=K>k_1$ such that homogeneous counterpart to (6.5) is asymptotically stable, other conditions of this statement hold, and $\limsup_{t\to\infty}\left\|V_K\left(t\right)\right\|=\bar{V}_K<\infty$; then,

$$\limsup_{t\to\infty}\left\|x\left(t,x_0\right)\right\|\le F_0\bar{V}_KZ_{K,\varphi,s}\left(t_0,\hat{z}_K\right),\ x_0\mid\left\|V_K^{-1}\left(t_0\right)x_0\right\|<\underline{z}_{K,\varphi},\ \hat{z}_K\in\left(0,\hat{\underline{z}}_{K,1}\right) \tag{6.14}$$

**Proof**. The proof of this statement can be attained directly by the utility of (6.6) and the apparent adjustment of the proof of Theorem 3 □

## 7. Simulations

This section applies the developed methodology to evaluate the degree of stability/ boundedness and estimate the evolutions of the norms of solutions to two nonlinear systems with time–dependent nonperiodic coefficients which are frequently appeared in various applications [32]. These systems comprise of two coupled Van der – Pol – like or Duffing - like oscillators with variable coefficients. The simulations authenticate our theoretical inferences and contrast the efficacy of the current and prior techniques [3] in their common application domains.

### 7.1 Coupled Van der Pol – like System

Equation (2.1) can be turned into such system if we assume that $B(t)=A+\mathrm{X}\left(t\right)$, where

$$A=\begin{pmatrix}0 & 1 & 0 & 0\\ -\left(\omega_1^2+\kappa\right) & -\alpha_1 & \kappa & 0\\ 0 & 0 & 0 & 1\\ \kappa & 0 & -\left(\omega_2^2+\kappa\right) & -\alpha_2\end{pmatrix},\ \mathrm{X}(t)=\begin{pmatrix}\begin{matrix}0 & 0\\ -\chi_{11} & -\chi_{12}\end{matrix} & 0\\ 0 & \begin{matrix}0 & 0\\ -\chi_{21} & -\chi_{22}\end{matrix}\end{pmatrix} \quad (7.1)$$

$$f_*=\begin{pmatrix}0 & -\varepsilon_1 x_2^3 & 0 & -\varepsilon_2 x_4^3\end{pmatrix}^T \quad (7.2)$$

$\alpha_1=0.4,\ \alpha_2=0.2,\ \omega_1^2=1,\ \omega_2^2=4$ , $\chi_{11}=\chi_{22}$, $\chi_{12}=\chi_{21}$, $\chi_{11}=a_{11}\sin\omega_{11}t+a_{12}\sin\omega_{12}t$ , $\chi_{12}=a_{21}\sin\omega_{21}t+a_{22}\sin\omega_{22}t$ , $F_{*i}=\bar{F}_i\sin q_i t,\ i=1,2$ , $q_1=5.43,\ q_2=10$ . In most but all our simulations we set that $\varepsilon_1=\varepsilon_2=0.5$ .

Note that the technique developed in [3] works under hypothesis that $A$ is a stable and $\mathrm{X}(t)$ is zero-mean matrices, whereas our current approach does not rest on these assumptions.

To simplify our task, we adopt in this section the appropriate formulas from Section 3 where we set that $k=1$.

In turn, let us recall that $v_1(t)=\left[v_{1,ij}(t)\right],\ i,j=1,\dots,4$ and $\Lambda_1(t)=diag\left(\lambda_1(t),\lambda_1^*(t),\lambda_2(t),\lambda_2^*(t)\right)$ are the eigenvectors and eigenvalues matrices of $a_0(t)$, $\lambda_k(t)=\alpha_k(t)+i\beta_k(t)$, $\lambda_k^*(t)=\alpha_k(t)-i\beta_k(t)$ , $k=1,2$, $G_1(t)=v_1^{-1}(t)v_1'(t)+v_1^{-1}(t)g_0(t)v_1(t), i,j=1,\dots,4$, $F_1(t)=v_1^{-1}(t)F_*(t)$ and

$$f_1=v_1^{-1}\left(0 \quad -\varepsilon_1\left(\sum_{k=1}^{4}v_{2k}y_k\right)^3 \quad 0 \quad -\varepsilon_2\left(\sum_{k=1}^{4}v_{4k}y_k\right)^3\right)^T.$$

Finally, we select $\bar{G}_1=G_1-\mathrm{Im}\left(diag\left(G_1\right)\right)$ and set that,

$$\|f_1\|\le L(t,z)=\left\|v_1^{-1}(t)\right\|\left(\varepsilon_1 e_1(t)+\varepsilon_2 e_2(t)\right)z^3 \quad (7.3)$$

where formulas for $e_i(t),\ i=1,2$ were derived in [3] under assumption that $v_1=const$; still these derivations remain intact if $v_1=v_1(t)$. Below we briefly replicate some of these steps for the sake of completeness.

Let $Y=\begin{pmatrix}0 & -\varepsilon_1 Y_2^3 & 0 & -\varepsilon_2 Y_4^3\end{pmatrix}^T$ , where $Y_i=\sum_{k=1}^{4}v_{ik}(t)y_k,\ i=2,4$ . Then,

$\|f_1\|=\left\|v_1^{-1}(t)Y\right\|\le\left\|v_1^{-1}(t)\right\|\|Y\|$ . Next, we adopt that $\|Y\|_2\le\|Y\|_1=abs\left(\varepsilon_1Y_2^3\right)+abs\left(\varepsilon_2Y_4^3\right)$, where $abs(u)=|u|$ if $\mathrm{Im}(u)=0$ and $abs(u)=\sqrt{a^2+b^2}$ if $u=a+ib$, and also that $\left\|y_k^m\right\|\le\|y\|^m$ , $m\in\mathbb{N}$ .

This yields that $\|Y\|_1=abs\left(\varepsilon_1\left(\sum_{k=1}^{4}v_{2k}(t)y_k\right)^3\right)+abs\left(\varepsilon_2\left(\sum_{k=1}^{4}v_{4k}(t)y_k\right)^3\right)=\left(\varepsilon_1e_1(t)+\varepsilon_2e_2(t)\right)\|y\|^3$,

where $e_1=\varepsilon_1\left(\sum_{k=1}^{4}abs\left(v_{2k}(t)\right)\right)^3$ , $e_2=\varepsilon_2\left(\sum_{k=1}^{4}abs\left(v_{4k}(t)\right)\right)^3$ which, finally, returns (7.3).

Afterwards, let us assess in simulations the efficacy of our optimal criterion which yields separation of temporal scales of matrix $B(t)$ via maximizing the degrees of stability of the linearized auxiliary equation (3.4) with $k=1$. This also let us to contrast the methodologies developed in [3] and in the current paper since the former – relays on averaging of $B(t)$ on infinite interval $[t_0 \quad \infty)$.

To reduce the scope of simulations, we set that $\delta_i = \bar{\delta},\ i = 1,...,4$ which embraces a suboptimal separation of temporal scales in $B(t)$. Figures 1a-f examine dependence of the Lyapunov exponents $\phi_{11}$ and $\phi_{12}$ upon $\bar{\delta}$, see Section 3 of this paper. For this sake, we used the following sets of parameters: in Figure 1a, $a_{11} = a_{12} = 0,$ $a_{21} = a_{22} = 0.5$ , $\omega_{21} = 20, \omega_{22} = 0.1$ ; in Figure 1b, $a_{11} = a_{12} = 0, a_{21} = a_{22} = 0.5, \omega_{21} = 20, \omega_{22} = 1$ ; in Figure 1c, $a_{11} = a_{12} = 0, a_{21} = a_{22} = 0.5, \omega_{21} = 20, \omega_{22} = 10$; in Figure 1d, $a_{11} = a_{12} = 0.5,\ a_{21} = a_{22} = 0$, $\omega_{11} = 20, \omega_{12} = 0.1$; in Figure 1e, $a_{11} = a_{12} = 0.5,\ a_{21} = a_{22} = 0,\ \omega_{11} = 20, \omega_{12} = 1$ and, lastly, in Figure 1f, $a_{11} = a_{12} = a_{21} = a_{22} = 0.5,\ \omega_{11} = 20, \omega_{12} = \omega_{22} = 1, \omega_{21} = 10$ .

Clearly, both Lyapunov exponents depend alike upon $\bar{\delta}$ and, as is expected, $\phi_{11}(\bar{\delta}) > \phi_{12}(\bar{\delta})$ on all figures. Still, the shapes of the plotted curves are significantly influenced by the values of the core frequencies of the variable coefficients used in simulations. The curves in Figures 1a and 1d as well as 1b and 1e are fairly similar. Figures 1a and 1d display that the current approach is somewhat superior then the one developed in [3] which averages $B(t)$ over interval $[t_0, \infty)$. Figure 1c shows that the technique developed in [3] are practically optimal for the set of parameters used in these simulations. Yet, Figures 1b, 1e, and 1f display that the current methodology frequently selects the optimal moving average with a window of finite length, which is determined by the core frequencies of the variable coefficients included in the relevant simulations.

Figures 2a -d displays in solid and dashed lines time-histories of slow, i.e., $\alpha_{1,\max}(t, \bar{\delta})$, and fast, $\left\| \bar{G}_1(t, \bar{\delta}) \right\|$, time-varying components of the linear part of equation (3.4) which are simulated for optimal value of $\bar{\delta}$ and the sets of parameters used in Figures 1a-c and 1f, respectively. Clearly, the proposed criterion provides a quite effective separation of the temporal scales of the linear block of the auxiliary equation despite that we evoke here its the most restricted and, hence, suboptimal version.

Note that all subsequent simulations of the auxiliary equations presented in this section include suboptimal values of $\bar{\delta}$ that were preliminary defined for the relevant sets of parameters.

Figure 3a displays in solid, dashed, and dashed-dotted lines the time histories of the norms of solutions to equation (2.2) and its upper bounds that are obtained in simulations of equations (3.4) and (3.2) with $F_1 = 0$ , respectively. These simulations assumed the following set of parameters, $\varepsilon_1 = \varepsilon_2 = 0.5,\ a_{11} = a_{12} = a_{21} = a_{22} = 0.5$, $\omega_{11} = 20,\ \omega_{12} = 1,\ \omega_{21} = 10, \omega_{22} = 0.5$. As is expected, solutions to equation (3.4) provide more efficient estimation of the norms of solutions to equation (2.2) than its more conservative counterpart. Figure 3b displays in solid and dashed lines the time histories of the norms of the solutions to equations (2.2) and (3.4) for the following set of parameters, $\varepsilon_1 = \varepsilon_2 = 0.5,\ a_{11} = a_{12} = 0,\ a_{21} = a_{22} = 0.1$ , $\omega_{21} = 10,\ \omega_{22} = 0.5$. Clearly, the reduction in the amplitudes of variable coefficients improves the accuracy of our estimates.

Figures 3c and 3d contrast in the solid and dashed lines the time histories of the norms of solutions to the nonhomogeneous equation (2.1) and its scalar counterpart (3.3), where we set that $\bar{F}_2 = 0$ and either $\bar{F}_1 = 0.025$ or $\bar{F}_1 = 0.05$, respectively. Note that the other parameters in these two figures and Figure 3b are the same. Clearly, the precisions of our estimates are correlated inversely with $\left\| \bar{F} \right\|$ .

Figure 3e displays in dashed and solid lines projections on $x_1 \times x_2$ -plane of the boundary of the stability region of the trivial solution to (2.2) that is developed in simulations of this equation and (3.4), respectively. Let us recall that simulations of the scalar equation (3.4) estimate the threshold value of $z_1(t, z_{10})$, i.e., $\bar{z}_1$ which is used, subsequently, to estimate the boundary of the stability region via application of the formula $\bar{z}_1 = \left\| V^{-1} x_0 \right\|$. The plots on this figure were developed for the followings set of parameters,
$\varepsilon_1 = \varepsilon_2 = 0.5,\ a_{11} = a_{12} = a_{21} = a_{22} = 0.1,\ \ \omega_1 = 20,\ \omega_{12} = 1, \omega_{21} = 10,\ \omega_{22} = 0.5$.

Figure 3f displays in dashed and solid lines the projections on $x_1 \times x_2$ -plane of boundaries of the trapping region which are developed in the simulations of equations (2.1) and (3.3), respectively. In this case we set that $\bar{F}_1 = 0.025, \bar{F}_2 = 0$. As prior, the latter simulations are used to estimate the threshold value of $z_1\left(t, z_{1,0}\right)$, i.e., $\bar{z}_1$. Clearly, the developed technique provides rather conservatively estimates of the boundary of the actual trapping region which, nonetheless, authenticate our theoretical inferences. Note that additional simulations provide alike estimates for projections on complementary coordinate frames.

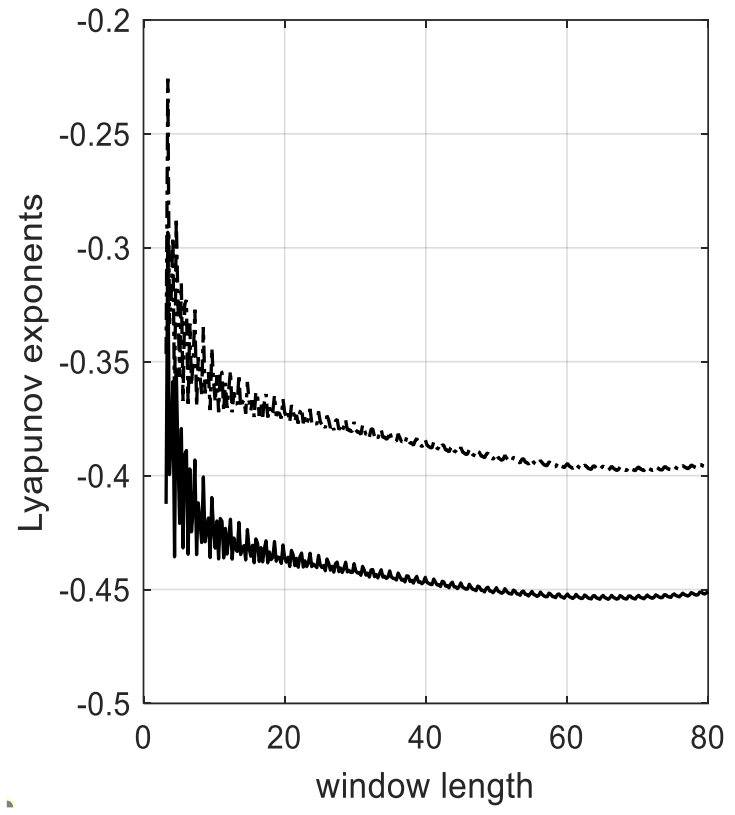

Fig. 1.a

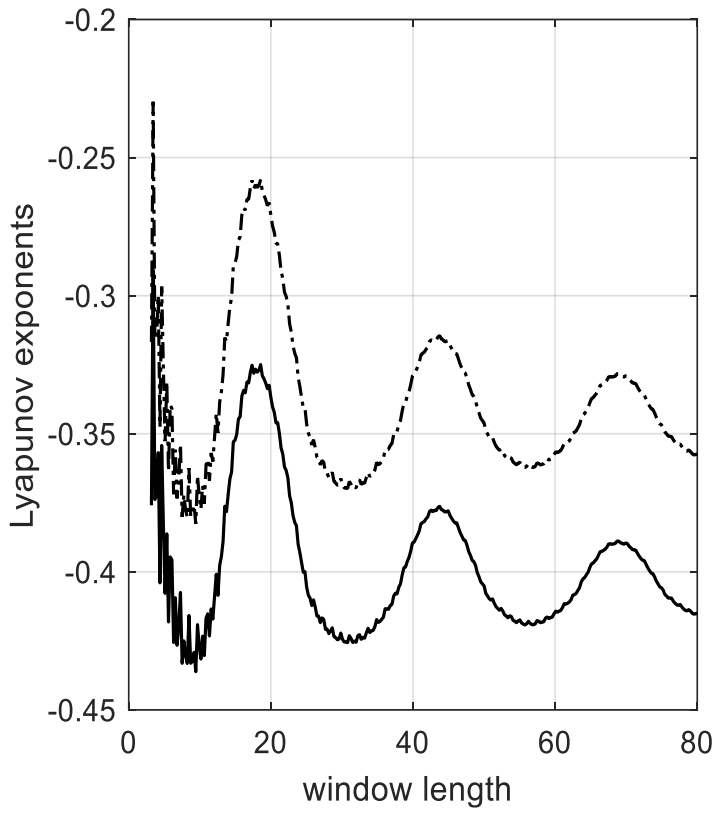

Fig. 1.b

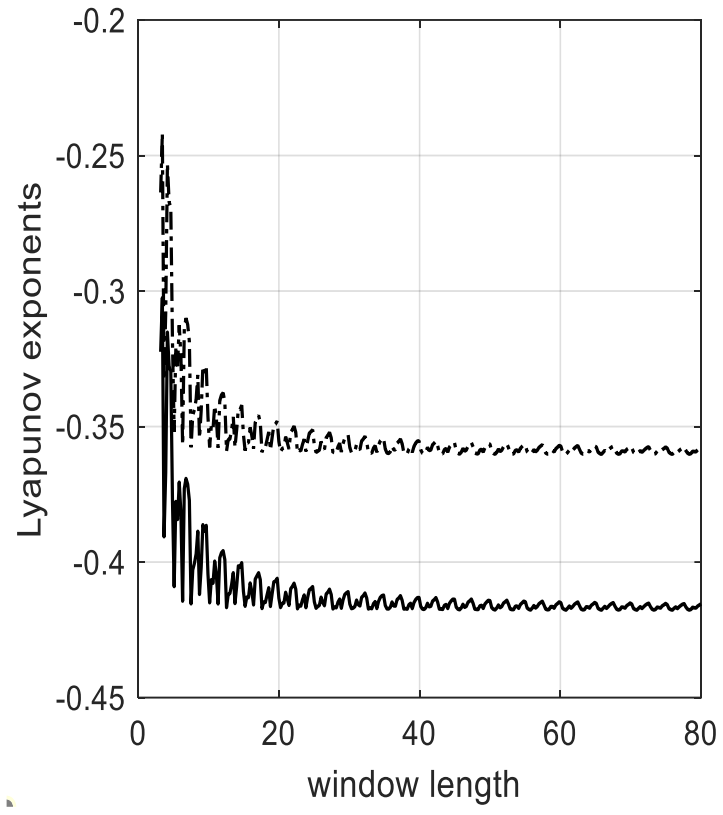

Fig. 1.c

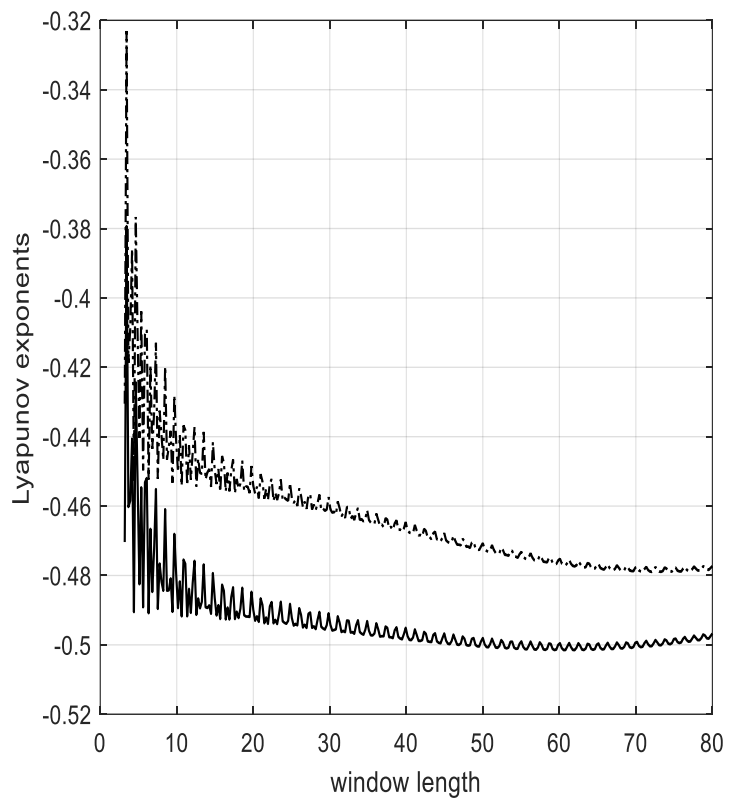

Fig. 1.d

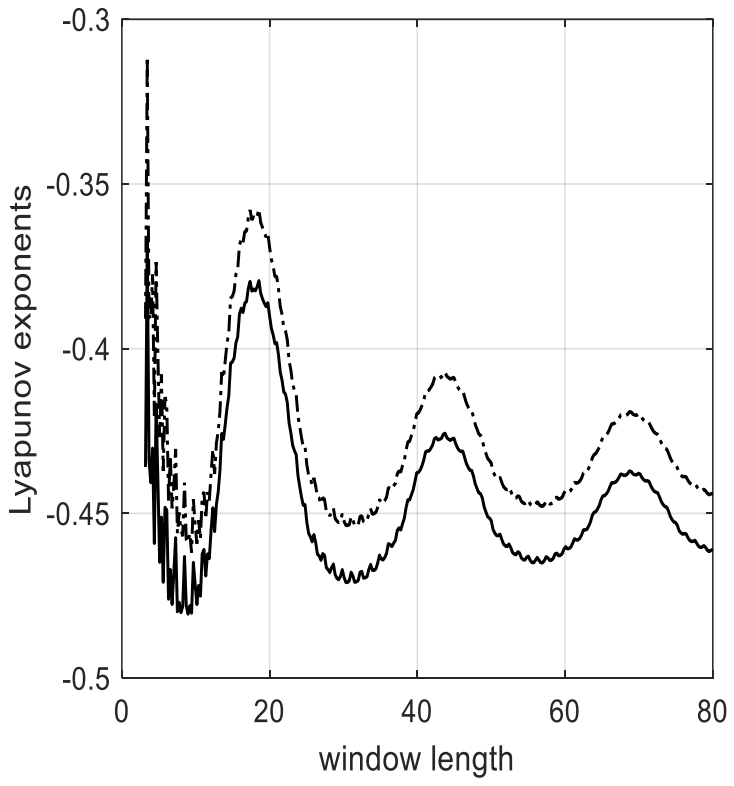

Fig. 1.e

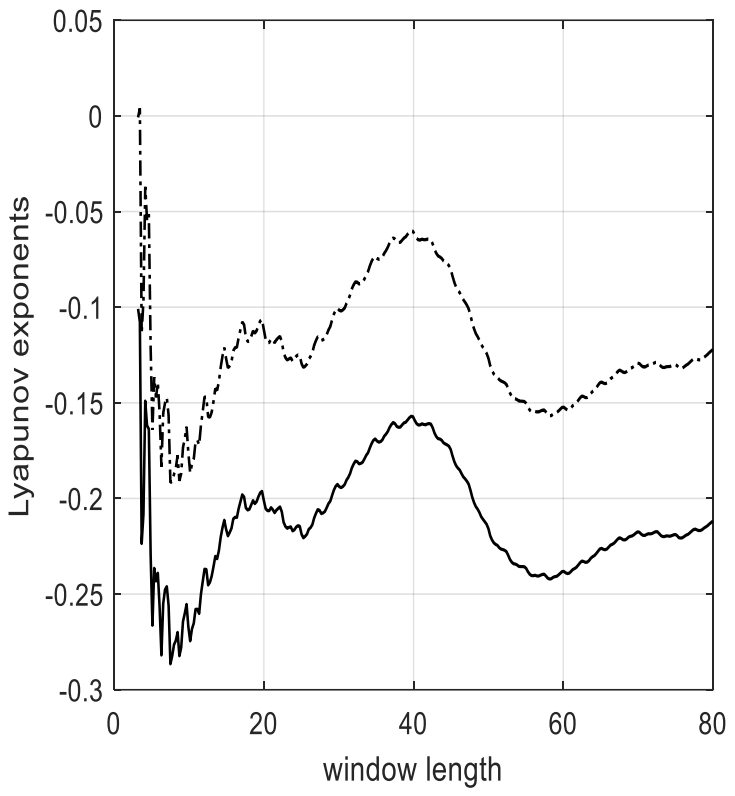

Fig. 1.f

Figures 1a-f plot in solid and dashed lines dependence of the Lyapunov exponents $\phi_{11}\left(\bar{\delta}\right)$ and $\phi_{12}\left(\bar{\delta}\right)$ upon $\bar{\delta}$ - the length of window of the moving average that is used for the separation of fast and slow time-varying components in entries of matrix $B(t)$.

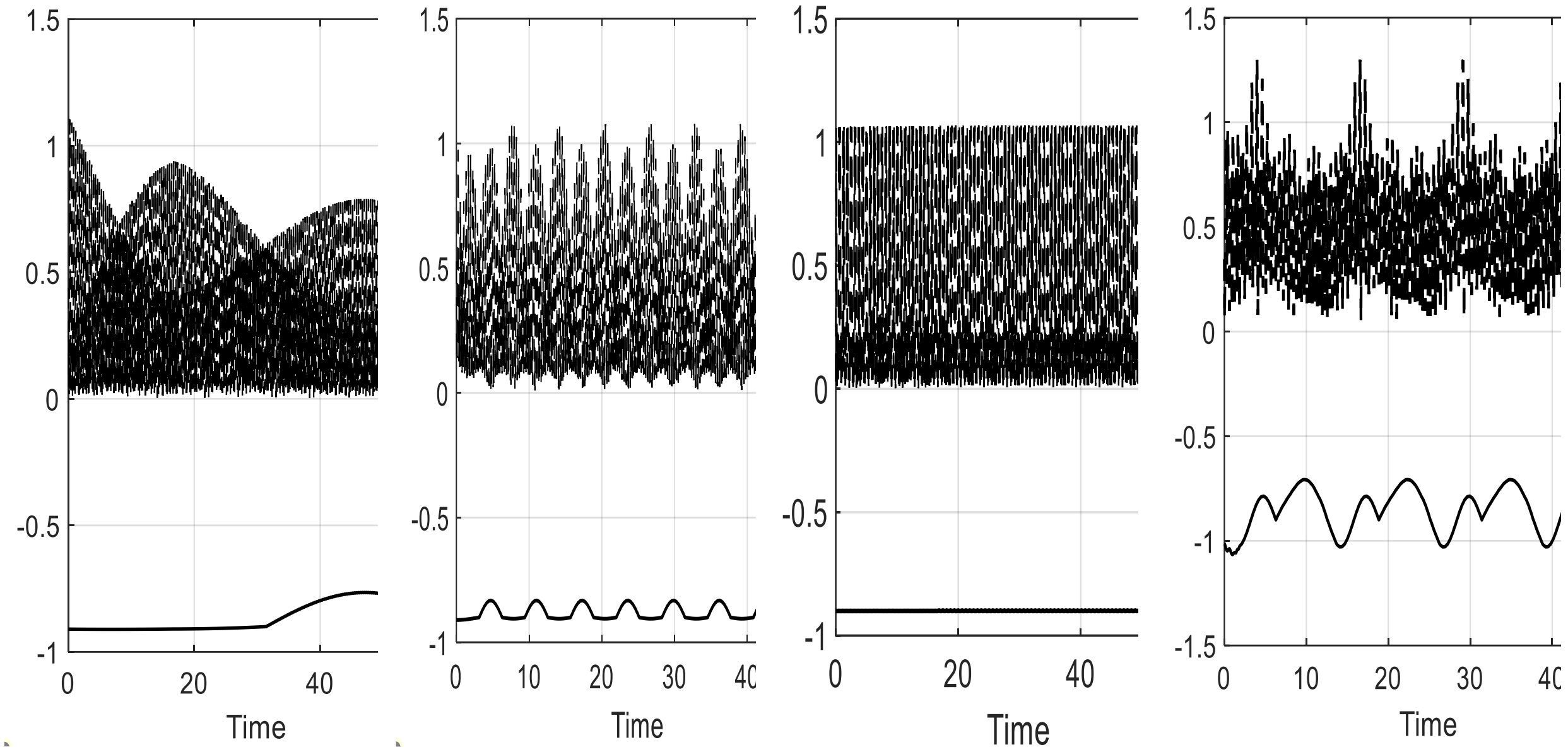

Fig. 2a Fig. 2b Fig. 2c Fig. 2d

Figures 2a -d display separation of the temporal scales attained in the linear block of equation (3.4) for the optimal value of $\bar{\delta}$ and the sets of parameters used in Figures 1a-c and 1f. Time-histories of slow, $\alpha_{1,\max}\left(t,\bar{\delta}\right)$, and fast, $\left\|\bar{G}_1\left(t,\bar{\delta}\right)\right\|$, components of the linear block of (3.4) are plotted in solid and dashed lines respectively.

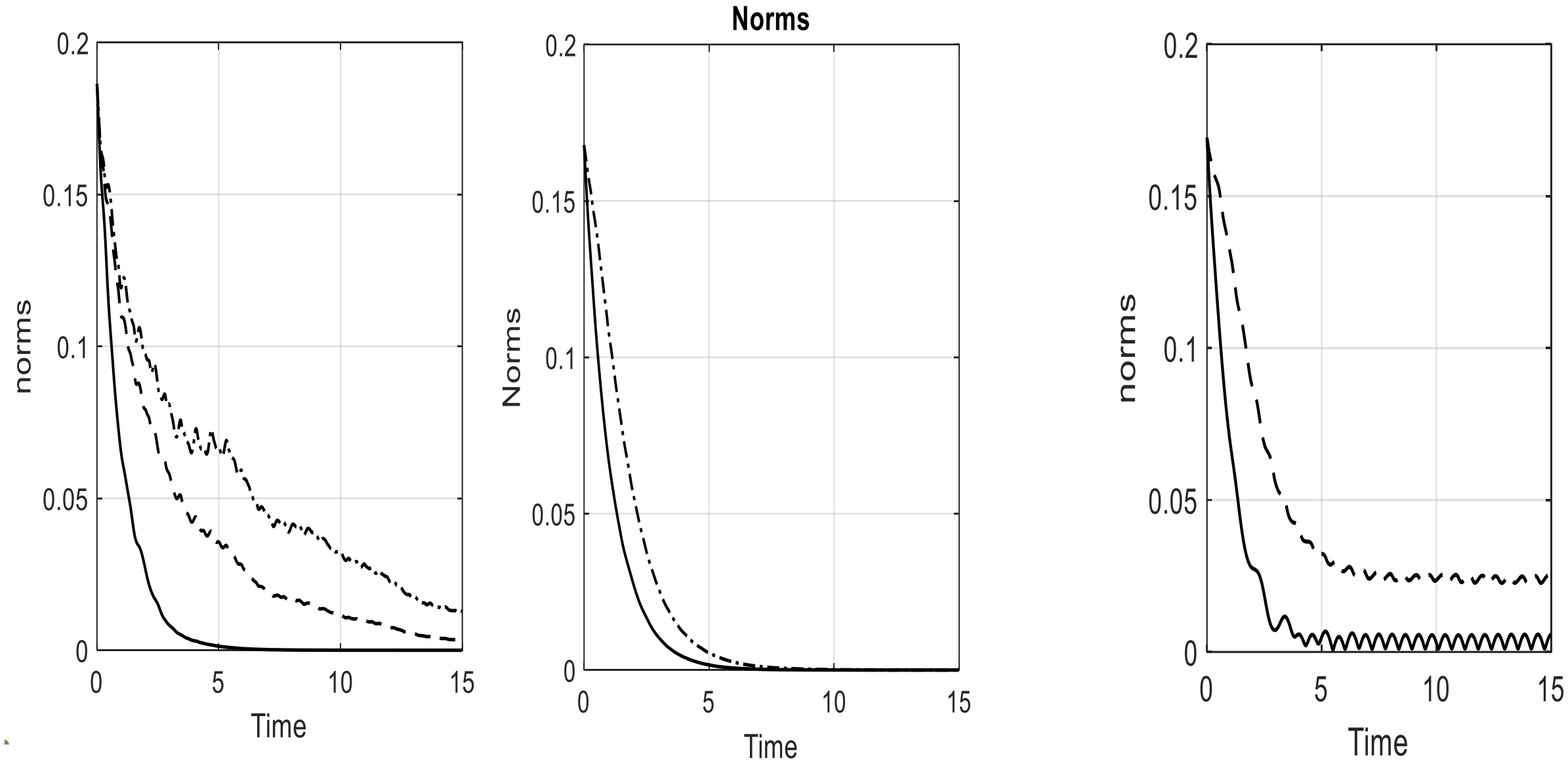

Fig. 3a Fig. 3b Fig. 3c

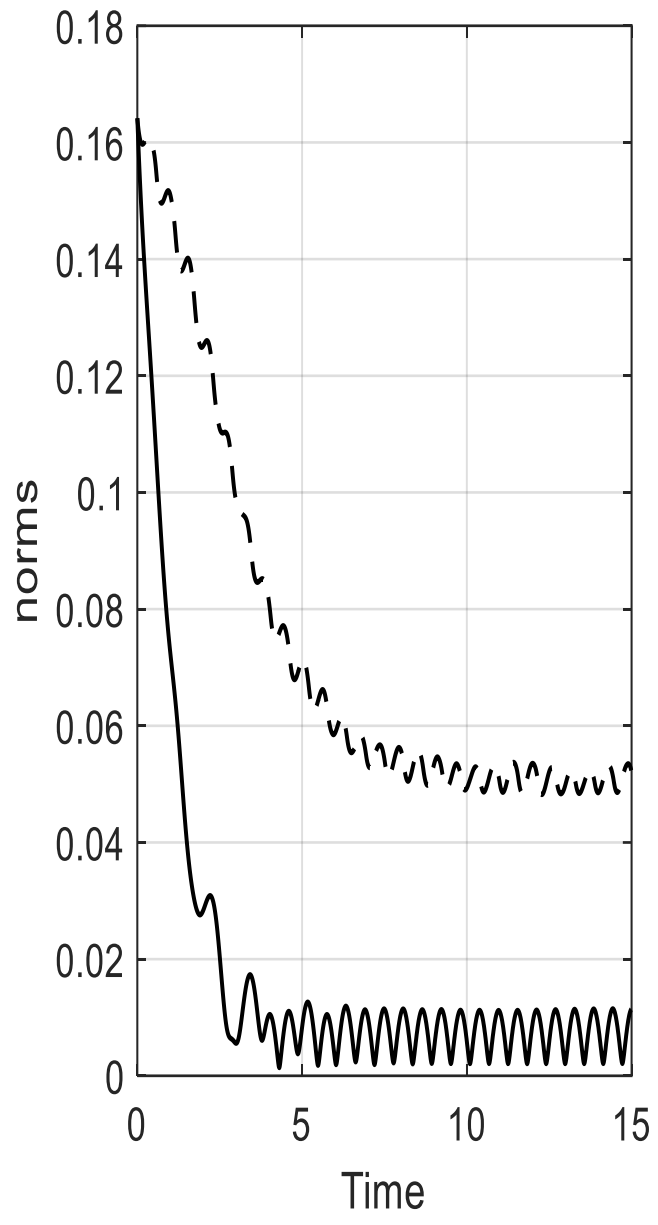

Fig. 3d

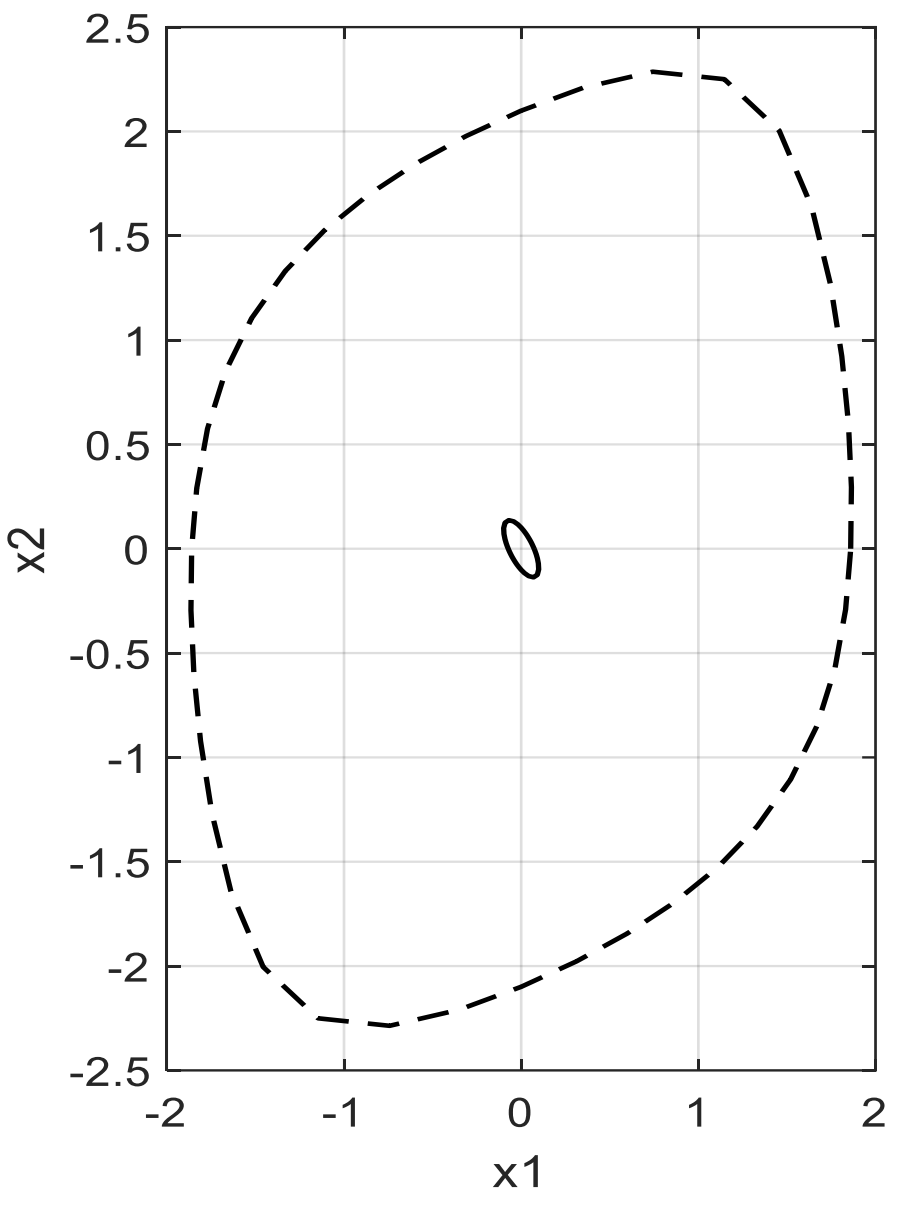

Fig. 3e

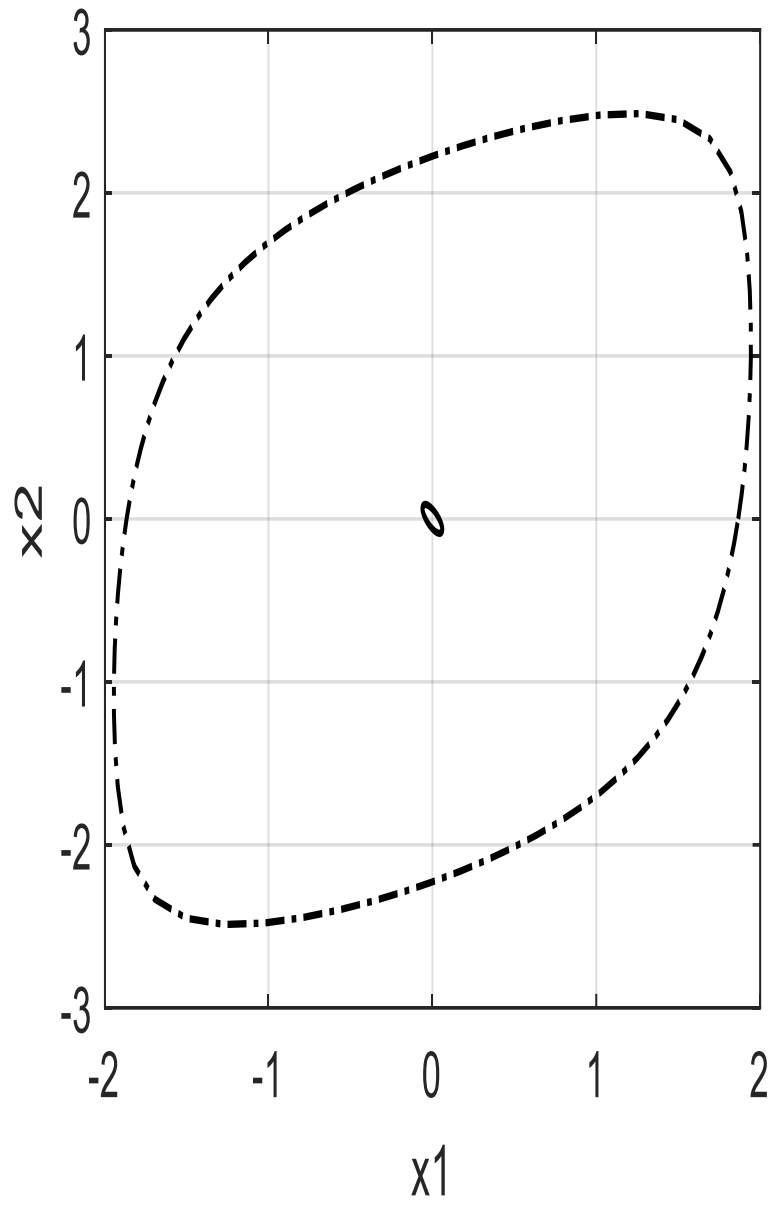

Fig. 3f

**7.2 Coupled Duffing-like System**

In a coupled Duffing-like system, the linear components are still defined by (7.1), whereas the definition of nonlinear components should be adjusted as follows,

$$f_* = \begin{pmatrix} 0 & -\varepsilon_1 x_1^3 & 0 & -\varepsilon_2 x_3^3 \end{pmatrix}^T$$

which prompts the relevant correction in (7.3), i.e., $v_{2k} \to v_{1k}$ and $v_{4k} \to v_{3k}$ . Figures 4a-f display the results of simulations of a coupled Duffing-like system and its scalar counterpart. Note that the curves in similar plots in Figures 3 and 4 are labeled alike.

In fact, Figures 4a contrasts in solid, dashed, and dotted dashed lines time -histories of the norms of solutions to equation (2.2) for the coupled Duffing-like system and their upper bounds which are developed in simulations of equations (3.4) and (3.2) with $F_1 = 0$ , respectively. As for the Van der Pol-like system, the solutions to equation (3.4) provide more efficient estimates than ones defined by (3.2) with $F_1 = 0$ . These simulations assumed the following set of parameters $a_{ij} = 0.5, i, j = 1, 2,\ \omega_{11} = 0.94,\ \omega_{12} = 2.12,\ \omega_{21} = 20$ , $\omega_{22} = 0.3,\ \varepsilon_1 = \varepsilon_2 = 0.5$.

Figure 4b displays in solid and dashed lines time histories of the norms of the solutions to equation (2.2) and its upper estimates which are developed in simulations of (3.4) with matched initial values for the following set of parameters, $\varepsilon_1 = \varepsilon_2 = 0.5,\ a_{11} = a_{12} = 0,\ a_{21} = a_{22} = 0.1$ , $\omega_{21} = 10,\ \omega_{22} = 0.5$ . As in Figure 3b, the amplitudes of the variable coefficients are scaled down in this figure. Nonetheless, the accuracy of our estimates improves marginally in this case.

Figure 4c contrasts in dashed and solid lines projections on $x_1 \times x_2$ -plane of the boundaries of the stability region of the trivial solution to (2.2) which were developed in simulations of this equation and its scalar counterpart (3.4), respectively. Note that in figures 4c -f, we set $a_{ij} = 0.1,\ i, j = 1, 2$ whereas other parameters adopted in these simulations are the same that the ones used in Figure 4a.

Furthermore, Figures 4d and 4e contrast in solid and dashed line time – histories of the norms of solutions to (2.1) and its scalar counterpart (3.3). In these figures we set that $\bar{F}_1 = 0.01,\ \bar{F}_2 = 0$ and $\bar{F}_1 = 0.025,\ \bar{F}_2 = 0$ , respectively. As is expected, enlargement of $\|\bar{F}\|$ decreases the accuracy of our estimates.

Figure 4f contrasts in dashed and solid lines projections on $x_1 \times x_2$ -plane of the boundaries of the trapping region of (2.1) which are developed in simulations of this equation and its scalar counterpart (3.3), respectively. In this case $\bar{F}_1 = 0.01, \; \bar{F}_2 = 0$ .

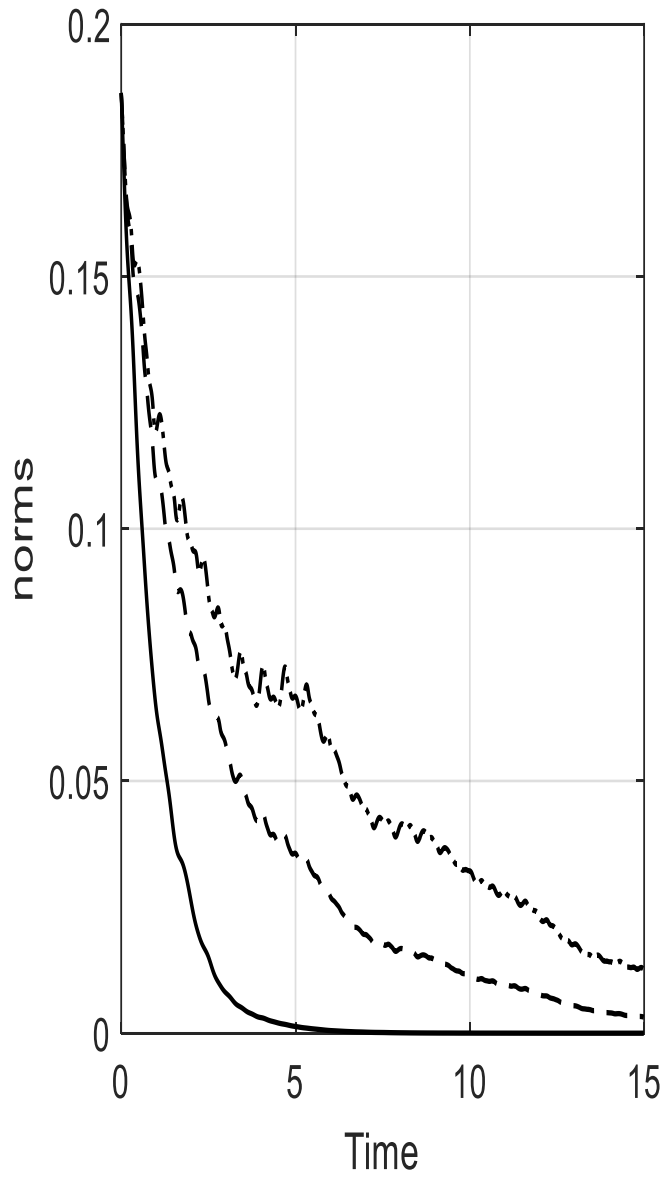

Fig. 4a

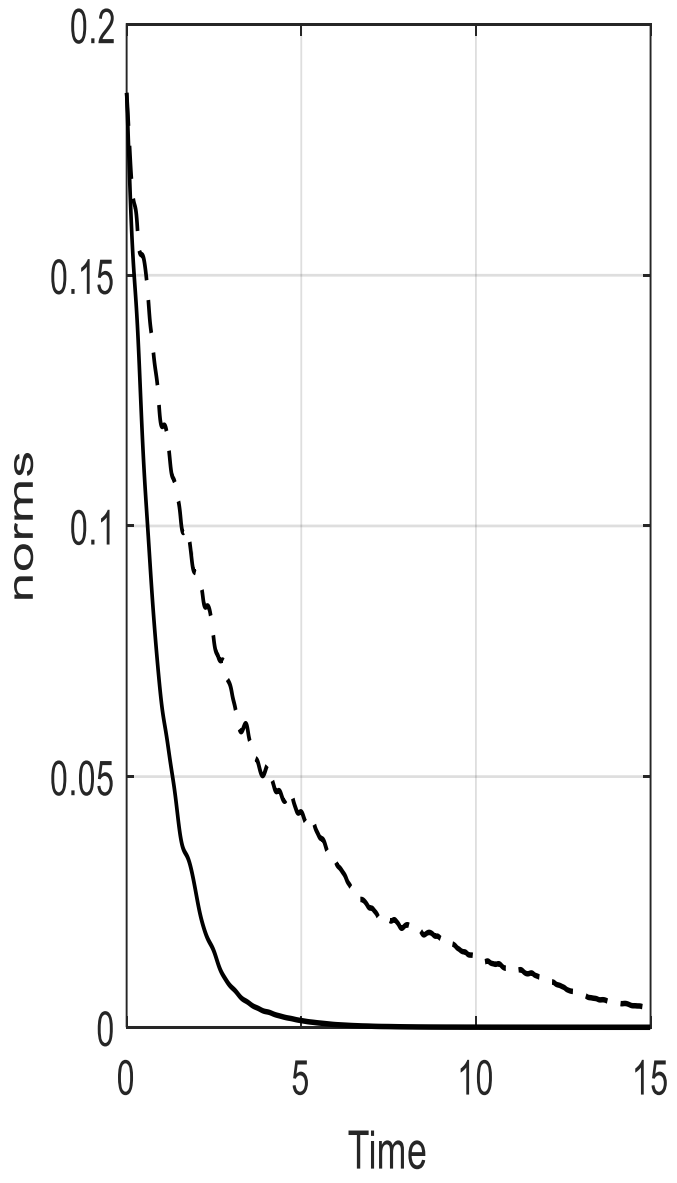

Fig. 4b

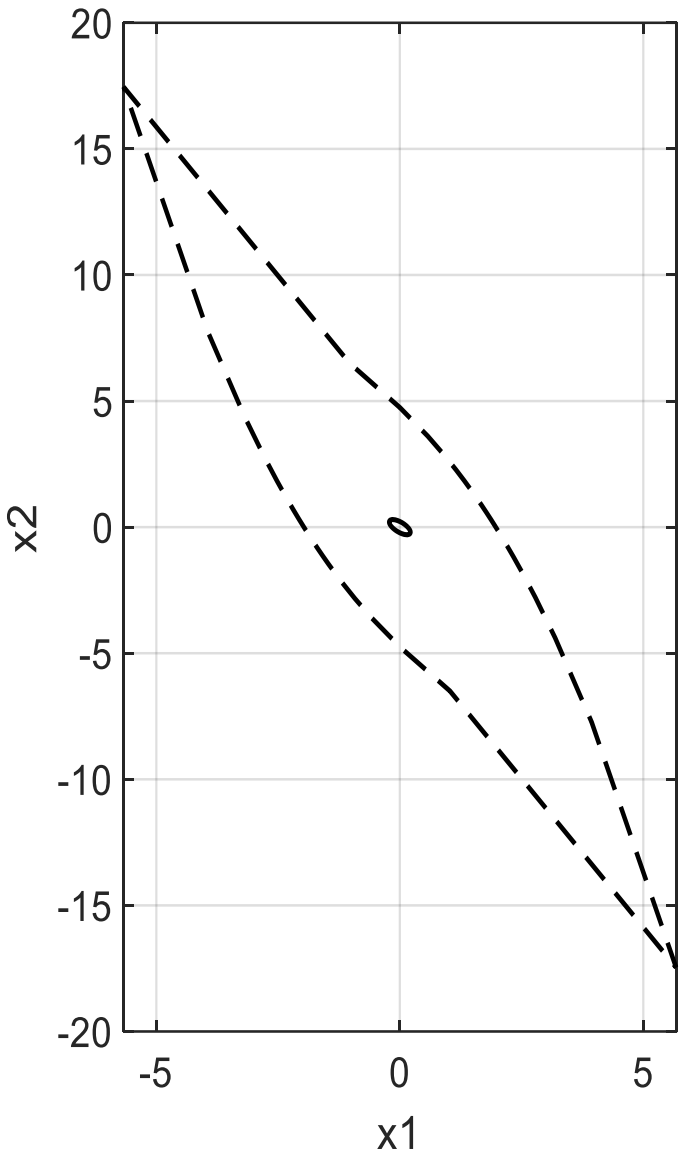

Fig. 4c

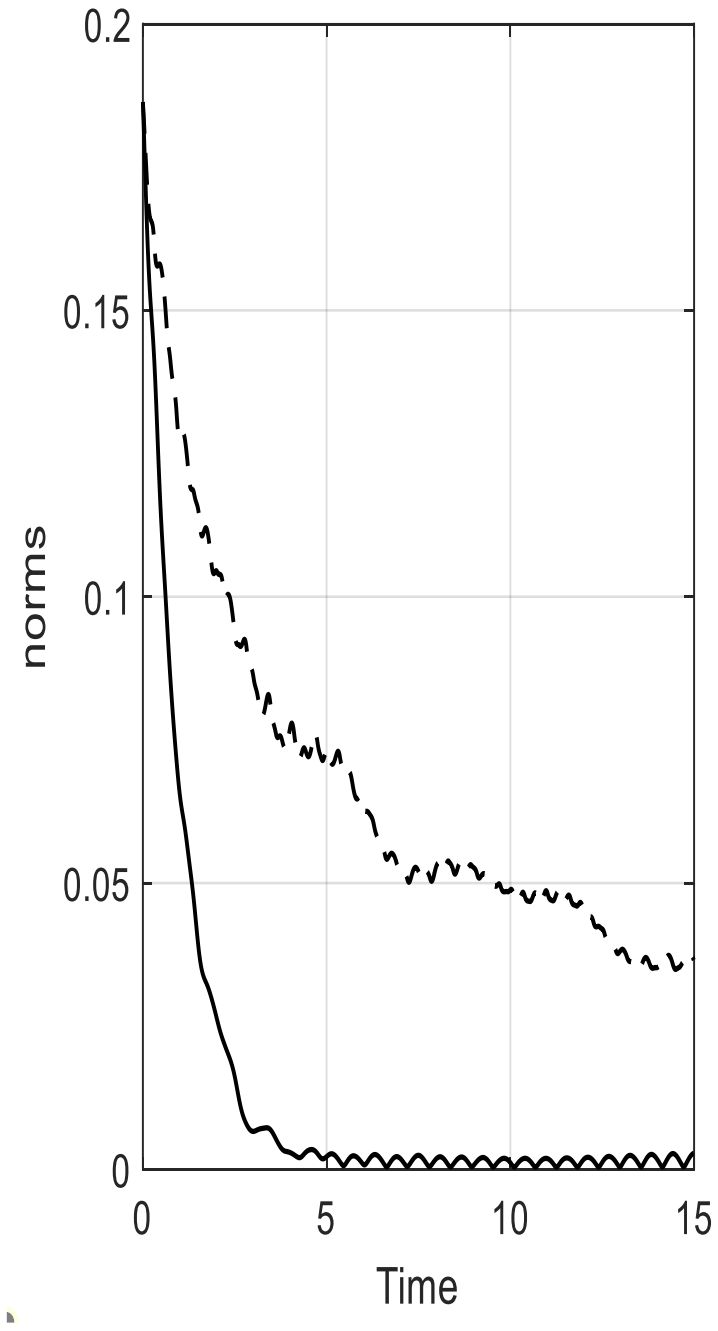

Fig. 4d

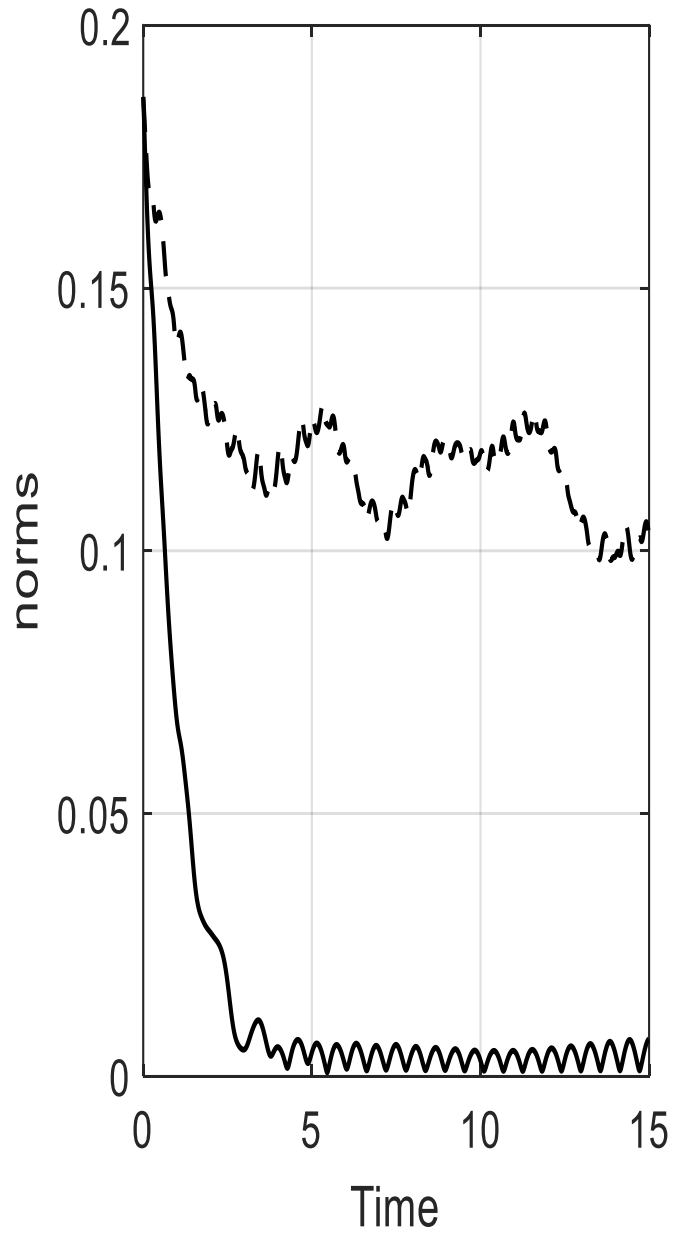

Fig. 4e

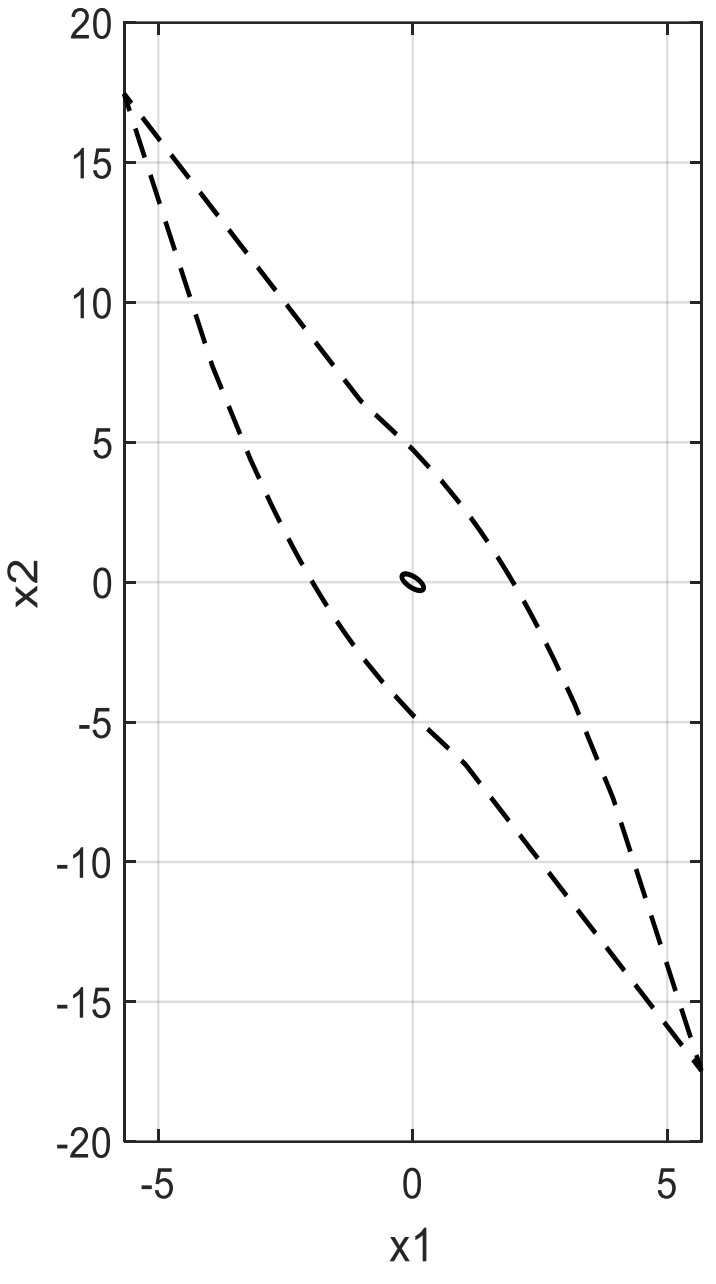

Fig. 4f

Apparently, the developed estimates of the norms of solutions to the coupled Van der-Pol-like and Duffing-like systems provide adequate accuracy if the solutions are emanated from a central part of the stability/ trapping regions. Yet, the boundaries of these regions are estimated rather conservatively by our current methodology. However, our estimates become more plausible if the original systems include some unknown but bounded in norm

components – a standard presumption in control theory. Still, all the above simulations affirmed our theoretical inferences.

## 8. Conclusions and Future Studies

This paper develops a novel methodology for estimation of the degree of boundedness/stability of nonlinear systems with time -dependent nonperiodic coefficients, which builds upon our recent studies in this area [1] - [3]. In [1], we develop a scalar auxiliary equation with solutions bounding from the above the norm of solutions to multidimensional nonlinear equations with time – varying coefficients, which let us to enhance the boundedness/stability criteria and estimate the corresponding regions for the underlying systems. However, some underlying assumptions constrain the broad application of this methodology, which, consequently, was recast in [3] under condition that the average of time -dependent matrix of the linear block of our system is a stable matrix. In contrast, the current approach is substantially more flexible. It merely requires splitting the time – dependent matrix of the linear block into slow and fast time-varying components. Such partitioning was developed by the application of moving averages to every entry of the underlying matrix. Next, the subsystem with a slow matrix is simplified using the relevant Lyapunov transform. Repeated applications of such transformations yield a sequence of abridged systems which leads to the development of their scalar auxiliary counterparts. Ultimately, the windows of moving averages are defined by maximizing the degree of stability of the linearized homogeneous auxiliary equation. This minimizes the conservatism of our estimates and aids the estimation of the behavior of solutions to the original systems through the assessment of the matching solutions to scalar auxiliary equations. Using this reasoning, we developed novel boundedness/stability criteria and estimated the boundedness/stability regions of the original systems.

Simplified estimates of solutions to scalar nonlinear auxiliary equations were developed using linear upper bounds for these equations. This let us to assess in closed form the evolutions of the norms of solutions to the original systems stemming from their trapping/stability regions and estimate the boundaries of these regions. The inferences developed by application of this technique should be useful in appropriate control tasks where some of them could reminiscent some of the estimates developed in the area of input-to-state stability under more restrictive conditions [33], [34].

The current methodology is simplified if the underlying system possesses only slow-varying coefficients. In this case, it can be used to authenticate and extend of practically plausible frozen coefficients or gain scheduling techniques which were frequently used in applications without appropriate justification.

Our inferences were examined in inclusive simulations that were partly discussed in this paper. These simulations contrast our current and prior methodologies and show that solutions to our scalar auxiliary equations deliver sound upper bounds of the norms of solutions to the original systems if they are stemmed from the central parts of their trapping/stability regions. However, the efficacy of our estimates is declined if the solutions are originated near the boundaries of these regions, which, in turn, are estimated rather conservatively. Nonetheless, all simulations reinforce our theoretical inferences and turn out to be more plausible for systems under uncertainty. Yet, integration of the methodologies developed in this paper and in [2] should substantially enhance estimation of the boundaries of the trapping/stability regions as well.

Acknowledgement. The programs used to simulate the models in Section 7 of this paper were developed by Steve Koblik.

Data sharing is not applicable to this article as no datasets were generated or analyzed during the current study.

The author declares that he has no conflict of interest.